\documentclass[twoside,11pt]{article}

\usepackage[margin=1in]{geometry}
\usepackage{mathptmx} 
\usepackage{authblk} 
\usepackage{setspace} 
\usepackage{parskip} 
\usepackage{fancyhdr} 
\usepackage{lastpage} 
\usepackage{amsmath,amssymb,amsthm}
\usepackage{bm} 
\usepackage{mathtools} 
\theoremstyle{plain}
\newtheorem{theorem}{Theorem}[section]
\newtheorem{lemma}[theorem]{Lemma}
\newtheorem{proposition}[theorem]{Proposition}

\theoremstyle{definition}
\newtheorem{definition}[theorem]{Definition}

\usepackage{graphicx}
\usepackage{subcaption} 
\usepackage{booktabs} 
\usepackage{multirow} 
\usepackage{multicol} 
\usepackage{longtable} 
\usepackage{rotating} 
\usepackage{caption} 
\usepackage{float} 
\graphicspath{{figures/}} 
\usepackage[
  backend=biber,       
  sorting=nty,          
  maxnames=3,           
  minnames=1,           
  maxbibnames=5,        
  minbibnames=3,       
  doi=true,             
  isbn=false,           
  url=false,            
  eprint=false          
]{biblatex}
\DeclareFieldFormat{journaltitle}{\textit{#1}}  
\DeclareFieldFormat{booktitle}{\textit{#1}}    
\DeclareFieldFormat{title}{#1.}

\newcommand{\keywords}[1]{\noindent\textbf{Keywords:} #1}

\usepackage[
  colorlinks=true,
  linkcolor=blue,
  citecolor=blue,
  urlcolor=blue,
  pdfauthor={Yiyuan Wang},
  pdftitle={A Nested Iterative Algorithm for Zero-Sum Linear-Quadratic Stochastic Differential Games in an Infinite Horizon},
  pdfsubject={Computational Approach; Zero-Sum Linear-Quadratic Stochastic Differential Games},
  pdfkeywords={Computational Approach; Zero-Sum Linear-Quadratic Stochastic Differential Games}
]{hyperref}
\usepackage{doi} 
\usepackage{url}

\title{A Nested Iterative Algorithm for Zero-Sum Linear-Quadratic Stochastic Differential Games in an Infinite Horizon}
\author{Yiyuan Wang \thanks{Yiyuan Wang is at the Shandong University Zhongtai Securities Institute for Financial Studies, Shandong University, 27 Shanda Nanlu, Jinan, P.R. China, 250100 (Email: wangyiyuan@mail.sdu.edu.cn).}}

\date{}

\begin{document}
\maketitle

\begin{abstract}
This paper proposes a new algorithm to compute closed-loop saddle points for infinite-horizon zero-sum linear-quadratic stochastic differential games via structural decoupling. Specifically, we develop a nested iterative scheme that generates a monotonically increasing matrix sequence to decompose the original problem into coupled subproblems. By sequentially solving the stabilizing solutions of associated algebraic Riccati equations for each subproblem, we recover the original problem’s stabilizing solution and rigorously prove sequence convergence. A numerical example further validates the effectiveness of the proposed method. To the best of our knowledge, this work extends the classical setting and provides the first general-purpose computational approach for this class of problems.\\
\vspace{0.5em}
\keywords{Computational Approach, Stabilizing Solution, Zero-Sum Linear-Quadratic Stochastic Differential Games, Game-Theoretic Algebraic Riccati Equations} \\
\end{abstract}

\section{Introduction} 
\label{sec:introduction}

Infinite-horizon two-player zero-sum linear-quadratic stochastic differential games (ZSLQSDGs) constitute a fundamental mathematical framework for modeling adversarial decision-making under uncertainty in continuous-time dynamical systems, with broad applications in robust control, mathematical finance, network security, and multi-agent coordination; see \cite{Bacsar2008,Alpcan2010,Savku2022}.

The theory of differential games originates in the seminal works of \cite{Isaacs1965} and \cite{Pontryagin1967}. \cite{Fleming1989} later extended the deterministic framework of \cite{Evans1984} to stochastic settings. Within the linear-quadratic (LQ) regime, systematic studies were carried out in \cite{Bacsar1998,Engwerda2005,Mou2006,Sun2014,Yu2015} and the references therein.
A landmark result was established by \cite{Sun2016}, who showed that a closed-loop saddle point exists if and only if a sign-indefinite game-theoretic algebraic Riccati equation (GTARE) admits a stabilizing solution; see also \cite{Sun2020_book,Dragan2020,Sun2021,Li2021,Goldys2022,Yu2022,Zhou2022,Zhang2024,Sun2025} for further developments. Consequently, computing a saddle point reduces to finding the stabilizing solution of the GTARE. However, standard numerical methods for conventional algebraic Riccati equations (AREs), such as that of \cite{Kleinman1968}, cannot be directly applied due to the sign-indefiniteness of the quadratic term.

For deterministic GTAREs, \cite{Lanzon2008} proposed a recursive method for positive stabilizing solutions. However, obtaining stabilizing solutions becomes substantially more intricate for stochastic systems subject to state-dependent noise;
\cite{Feng2010} addressed this gap by developing a method for GTAREs in a specific class of ZSLQSDGs. This work was extended by \cite{Dragan2011} to stochastic $H_\infty$ problems with state-dependent and control-dependent noise, and further generalized in \cite{Dragan2015} to include disturbance-dependent noise.
Later, \cite{Ivanov2018} studied a specific class of GTAREs, explored three computational approaches, and compared their numerical performance. Additional related work can be found in \cite{Aberkane2023, Sun2023} and their references. 
Despite this extensive literature, no general-purpose computational approach exists for the full class of infinite-horizon ZSLQSDGs with state-dependent and control-dependent noise for both players—the exact problem studied in this paper. The coupling structure poses a fundamental challenge when the diffusion term of the state equation depends simultaneously on the controls of both players: in this case, the blocks $R_{ij}(P)$ in the GTARE inherit cross terms $D_{l,i}^\top P D_{l,j}$ for $i\neq j$, producing a non-trivial interaction that cannot be resolved by existing single-player or single-control-in-noise techniques.

The key observation underlying this paper is that the coupled two-player game admits a natural decoupling: at each iteration, fixing the linear feedback strategies of both players---as derived from the previous iterate rather than assumed optimal---transforms the original game into a single-player optimal control problem whose optimality condition is a standard ARE with a sign-definite quadratic term. By iterating over a carefully constructed sequence of such decoupled subproblems, we obtain a nested iterative scheme that generates a monotonically increasing matrix sequence converging to the stabilizing solution of the original GTARE. The method requires only repeated solutions of standard AREs.

The principal contributions of this paper are threefold:
\begin{itemize}
    \item \textbf{Conceptual.} We identify and exploit the decoupling structure inherent in ZSLQSDGs, showing that the sign-indefinite GTARE can be decomposed into a hierarchy of sign-definite AREs.
    \item \textbf{Algorithmic.} We propose the first nested iterative framework capable of computing closed-loop saddle points for the full class of infinite-horizon ZSLQSDGs with state-dependent and control-dependent noise for both players. The framework subsumes deterministic GTAREs \cite{Lanzon2008} and stochastic $H_\infty$-type GTAREs \cite{Feng2010,Dragan2011,Dragan2015} as special cases.
    \item \textbf{Theoretical.} We establish complete convergence guarantees: under standard stabilizability and detectability assumptions, every subproblem admits a unique  stabilizing solution, and the nested iterates converge monotonically to the unique stabilizing solution of the GTARE.
\end{itemize}

\noindent The remainder of this paper is organized as follows.
Section~\ref{sec:preliminaries} introduces notation, formulates the problem, presents the decoupling representation, and describes the iterative algorithm.
Section~\ref{sec:main_results} develops the structural properties of the game operator, establishes well-posedness criteria for the subproblem AREs, and proves the main convergence theorem.
Section~\ref{sec:example} provides a numerical illustration.

\section{Preliminaries} 
\label{sec:preliminaries}

\subsection{Notation}
Let us first introduce the following notation:
\begin{itemize}
    \item $\mathbb{R}^n$: $n$-dimensional real Euclidean space; $\mathbb{C}^{-}$: the set of complex numbers with negative real part; $\mathbb{R}^{n\times m}$: the space of $n\times m$ real matrices; $\mathbb{S}^n$: the set of all $n\times n$ symmetric matrices; $\overline{\mathbb{S}}_{+}^n$: the set of all $n\times n$ symmetric positive semidefinite matrices; $\mathbb{S}_{+}^n$: the set of all $n\times n$ symmetric positive definite matrices. 
    \item $I_n$: the identity matrix of size $n$; $\mathbb{O}_{n \times m}$: the null matrix of size $n \times m$. It can be simplified as 0 when no ambiguity is generated; $A_{m \times n}^{\times (l)}$ denotes the object represented by $A_{m \times n}$ arranged repeatedly for $l$ times.
    \item $A^{\top}$: the transpose of the matrix $A$; $\langle\cdot,\cdot\rangle$: the inner product on a Hilbert space. If $A\in\mathbb{S}_{+}^n$ (resp., $A\in\overline{\mathbb{S}}_{+}^n$), we write $A\succ0$ (resp., $A\succeq0$). For any $A,B\in\mathbb{S}^n$, we use the notation $A\succ B$ (resp., $A\succeq B$) to indicate that $A - B\succ0$ (resp., $A - B\succeq0$).
    \item Let $\mathbb{H}$ be a Euclidean space, and we define the following space \footnote{$\varphi \in \mathbb{F}$ denotes that $\varphi$ is $\mathbb{F}$-progressively measurable.}: 
    $L_{\mathbb{F}}^2(\mathbb{H}) = \{\varphi: [0,\infty) \times \Omega \rightarrow \mathbb{H} \mid \varphi \in \mathbb{F},\mathbb{E} \int_{0}^{\infty} |\varphi(t)|^2 \, dt < \infty \} $; $\mathcal{U}_i = L_{\mathbb{F}}^2(\mathbb{R}^{m_i})(i = 1,2)$, $\mathcal{X}_{\mathrm{loc}}[0,\infty) = \{\varphi: [0,\infty) \times \Omega \rightarrow \mathbb{R}^n \mid \varphi \in \mathbb{F} \,\text{ is continuous},\, \mathbb{E}\left[\sup_{0 \leq t \leq T} |\varphi(t)|^2\right] < \infty, \forall\, T > 0 \}$; $\mathcal{X}[0,\infty) = \{\varphi \in \mathcal{X}_{\mathrm{loc}}[0,\infty) \mid \mathbb{E} \int_{0}^{\infty} |\varphi(t)|^2 \, dt < \infty \}$.
\end{itemize}

\subsection{Zero-Sum Linear-Quadratic Stochastic Differential Games in an Infinite Horizon}

Let $(\Omega,\mathcal{F},\mathbb{F},\mathbb{P})$ be a complete filtered probability space on which a r-dimensional standard Brownian motion $W = \{W(t)^{\top}=(w_{1}(t),\dotsb,w_{r}(t));t\geq0\}$ is defined with $\mathbb{F}=\{\mathcal{F}_t\}_{t\geq0}$ being the usual augmentation of the natural filtration generated by $W$. We consider the following controlled linear stochastic differential equation (SDE) on $[0,\infty)$:
\begin{equation} 
\label{zslqsdg:sde}
    \begin{cases} 
        &dX(t) = \lbrack AX(t) + B_1 u_1(t) + B_2 u_2(t) \rbrack dt \\
        &\quad\quad\quad+ \sum_{l=1}^{r}\lbrack C_{l}X(t) + D_{l,1} u_1(t) + D_{l,2} u_2(t) \rbrack dw_{l}(t) \\
        &X(0) = x
    \end{cases} ,
\end{equation}
in which $x \in\mathbb{R}^n$, $A,C_{l}\in\mathbb{R}^{n\times n}$ and $B_i,D_{l,i}\in\mathbb{R}^{n\times m_i}$ ($i = 1,2;1 \leq l \leq r$). In the above, the state process $X$ is an $n$-dimensional vector, player 1 $u_1$ and player 2 $u_2$ are an $m_1$-dimensional vector and an $m_2$-dimensional vector, respectively. 

In this paper, Player 1 and Player 2 share the same performance functional:
\begin{equation}
\label{zslqsdg:performance_functional}
    J(x;u_1,u_2)=\mathbb{E}\int_{0}^{\infty}\Bigg[\Bigg\langle\begin{pmatrix}
    Q & S_1^{\top} & S_2^{\top}\\
    S_1 & R_{11} & R_{12}\\
    S_2 & R_{21} & R_{22}
    \end{pmatrix}\begin{pmatrix}
    X(t)\\
    u_1(t)\\
    u_2(t)
    \end{pmatrix},\begin{pmatrix}
    X(t)\\
    u_1(t)\\
    u_2(t)
    \end{pmatrix}\Bigg\rangle
    \Bigg]dt,
\end{equation}
where the weighting coefficients satisfy
\[
    Q\in\mathbb{S}^n,\quad R_{21}^{\top}=R_{12}\in\mathbb{R}^{m_1\times m_2},\quad S_i\in\mathbb{R}^{m_i\times n}, R_{ii}\in\mathbb{S}^{m_i}(i = 1,2).
\]
For $(x,u_1,u_2)\in\mathbb{R}^n\times\mathcal{U}_1\times\mathcal{U}_2$, the solution $X(\cdot;x,u_1,u_2)$ to the SDE in \eqref{zslqsdg:sde} may only exist in $\mathcal{X}_{loc}[0,\infty)$, which renders $J(x;u_1,u_2)$ ill-defined. We define the set of admissible controls as $\mathcal{U}_{ad}(x)=\{ (u_1,u_2)\in\mathcal{U}_1 \times \mathcal{U}_2\mid X(\cdot;x,u_1,u_2)\in\mathcal{X}[0,\infty)\}$. A pair $(u_1,u_2)\in\mathcal{U}_{ad}(x)$ is called an admissible control pair for the initial state $x$, and the corresponding $X(\cdot;x,u_1,u_2)$ is referred to as the admissible state process. In this case, $J(x, u_1, u_2)$ is clearly well-defined.

In this zero-sum game, Player 1 (\textit{the maximizer}) chooses control $u_1$ to maximize \eqref{zslqsdg:performance_functional}, while Player 2 (\textit{the minimizer}) chooses $u_2$ to minimize the same function. The problem is to find an admissible control pair $(u_1^*,u_2^*)$ that both players can accept. We denote the above-mentioned problem as $({SDG})^{\,0}_{\infty}$ for short. For a description of the $({SDG})^{\,0}_{\infty}$ problem, refer to \cite{Sun2020_book} . More detailed information can be found therein. 

The $({SDG})^{\,0}_{\infty}$ problem corresponds to following stochastic GTARE:
\begin{equation}
\label{zslqsdg:gtare}
    \begin{cases}
        Q(P)-S(P)^{\top}R(P)^{-1}S(P)=0\\
        R_{11}(P)\prec 0,\quad R_{22}(P) \succ 0 
    \end{cases},
\end{equation}
where,
\begin{equation}
\label{zslqsdg:gtare_auxiliarymatrix}
\begin{aligned}
    &Q(P)=PA + A^{\top}P+\sum_{l=1}^{r}C_{l}^{\top}PC_{l} + Q,\\
    &S(P)=\begin{bmatrix}
    S_1(P) \\
    S_2(P)
    \end{bmatrix}=\begin{bmatrix}
    B_1^{\top}P+\sum_{l=1}^{r}D_{l,1}^{\top}PC_{l}+S_1 \\
    B_2^{\top}P+\sum_{l=1}^{r}D_{l,2}^{\top}PC_{l}+S_2\\
    \end{bmatrix},\\
    &R(P)=\begin{bmatrix}
    R_{11}(P) & R_{12}(P)\\
    R_{21}(P) & R_{22}(P)\\
    \end{bmatrix},\quad R_{ij}(P)=R_{ij}+\sum_{l=1}^{r}D_{l,i}^{\top}PD_{l,j}\,(i,j = 1,2).\\
\end{aligned}
\end{equation}

\begin{definition}
\label{zslqsde:mean_square}
The system 
\begin{equation*}
    \begin{cases} 
        dX(t) =  A X(t) dt + \sum_{l=1}^{r} C_{l}X(t) dw_{l}(t) \\
        X(0) = x
    \end{cases},
\end{equation*}
denoted as system $(A,C_{1},\dotsb,C_{r})$ is called mean-square stable such that the solution 
satisfies \[\lim_{t\to\infty} \mathbb{E}[X(t)^{\top}X(t)] = 0 \,\,\text{for every initial state}\,\, x \in \mathbb{R}^n.\]
\end{definition}

A concept related to the above admissible control pair is that of mean-square stabilizers. We define the set of all such stabilizers associated with the SDE \eqref{zslqsdg:sde} as
\[ 
\mathcal{K} = \left\{
    \begin{aligned}
        &\text{All of} \,\,\,(\varTheta_1, \varTheta_2) \in \mathbb{R}^{m_1 \times n} \times \mathbb{R}^{m_2 \times n} \,\text{such that the system} \,(A + B_1 \varTheta_1 + B_2 \varTheta_2,\\
        &C_1 +D_{1,1} \varTheta_1 + D_{1,2} \varTheta_2, \dotsb, C_r + D_{r,1} \varTheta_1 + D_{r,2} \varTheta_2) \text{ is mean-square stable}.
    \end{aligned}
\right\}.
\]
Let $(\varTheta_1, \varTheta_2) \in \mathcal{K}$, the feedback control $U(\cdot)=\begin{bmatrix} \varTheta_1 X(\cdot)\\ \varTheta_2 X(\cdot)\end{bmatrix}$ is called mean-square stabilizing. Moreover, for $\varTheta_i \in \mathbb{R}^{m_i \times n}(i = 1, 2)$, we let
\begin{align*}
&\mathcal{K}_1(\varTheta_2) = \left\{ \varTheta_1 \in \mathbb{R}^{m_1 \times n} : (\varTheta_1 , \varTheta_2 ) \in \mathcal{K} \right\}, \\
&\mathcal{K}_2(\varTheta_1) = \left\{ \varTheta_2 \in \mathbb{R}^{m_2 \times n} : (\varTheta_1 , \varTheta_2 )\in \mathcal{K} \right\}.\end{align*}
We refer to \cite{Rami2000, Sun2020_book} for the above definition.

\begin{definition}[\cite{Sun2020_book}]
A 2-tuple pair $(\varTheta_1^*, \varTheta_2^*) \in \mathbb{R}^{m_1 \times n} \times \mathbb{R}^{m_2 \times n}$ is called a closed-loop saddle point of Problem $({SDG})^{\,0}_{\infty}$ if $(\varTheta_1^*, \varTheta_2^*) \in \mathcal{K}$ and
\[J(x,\varTheta_1X(\cdot), \varTheta_2^*X(\cdot)) \leq J(x,\varTheta_1^*X(\cdot), \varTheta_2^*X(\cdot)) \leq J(x,\varTheta_1^*X(\cdot), \varTheta_2X(\cdot)),\]
for all $x \in \mathbb{R}^n, (\varTheta_1, \varTheta_2) \in \mathcal{K}_1(\varTheta_2^*) \times \mathcal{K}_2(\varTheta_1^*)$.
\end{definition}

\begin{definition}[\cite{Sun2020_book}]
A matrix $P\in\mathbb{S}^n$ is called a stabilizing solution of stochastic GTARE \eqref{zslqsdg:gtare} if $P$ satisfies \eqref{zslqsdg:gtare} and $(K_1(P), K_2(P)) \in \mathcal{K}$, where $\begin{bmatrix} K_1(P) \\ K_2(P) \end{bmatrix}=-R(P)^{-1}\begin{bmatrix}S_1(P) \\S_2(P)\end{bmatrix}$ and $\mathcal{K}$ is the set of all such stabilizers associated with the SDE \eqref{zslqsdg:sde}.
\end{definition}

To facilitate subsequent analysis of the problem, we introduce the following decoupling representation method. Setting 
\[
\begin{bmatrix}  v_1(t) \\ v_2(t) \end{bmatrix}=\begin{bmatrix}  u_1(t) \\ u_2(t) \end{bmatrix}-\begin{bmatrix} K_1(0) \\ K_2(0) \end{bmatrix},\quad 
\begin{bmatrix} K_1(0) \\ K_2(0) \end{bmatrix} = -R(0)^{-1} \begin{bmatrix} S_1(0) \\ S_2(0) \end{bmatrix},
\] and $v_2(t) = Lx(t)$ in $({SDG})^{\,0}_{\infty}$, we obtain
\begin{equation*} 
    \begin{cases} 
        dX(t) = \lbrack A_{L}X(t) + B_1v_1(t) \rbrack dt + \sum_{l=1}^{r}\lbrack C_{lL}X(t) + D_{l,1}v_1(t) \rbrack dw_{l}(t) \\
        X(0) = x
    \end{cases}
\end{equation*}
in which $x \in\mathbb{R}^n$, and 
\begin{equation*}
    J_{L}(x;v_1)=\mathbb{E}\int_{0}^{\infty}\Bigg[\Bigg\langle\begin{pmatrix}
    Q_{L} & S^{\top}_{L} \\
    S_{L} & R_{11}
    \end{pmatrix}\begin{pmatrix}
    X(t)\\
    v_1(t)
    \end{pmatrix},\begin{pmatrix}
    X(t)\\
    v_1(t)
    \end{pmatrix}\Bigg\rangle
    \Bigg]dt,
\end{equation*}
where
\[
\begin{cases}
&A_{L} = A + B_1K_1(0) + B_2K_2(0)+ B_2L\\
&C_{lL} = C_{l} + D_{l,1}K_1(0) + D_{l,2}K_2(0)+ D_{l,2}L,  1 \leq l \leq r \\
&Q_{L} = Q - S^{\top}(0) R(0)^{-1}S(0) + L^{\top}R_{22}L,\quad S_{L}= R_{12}L 
\end{cases}.
\]
The corresponding ARE is
\begin{equation}
\label{gtare_lw}
\begin{aligned}
&PA_{L} + A_{L}^{\top}P+ \sum_{l=1}^r C_{lL}^{\top}PC_{lL}+ Q_{L} - (B_{1}^{\top}P + \sum_{l=1}^r D_{l,1}^{\top}PC_{lL} + S_{L}) ^{\top} \\
& \quad \times(R_{11} + \sum_{l=1}^r D_{l,1}^{\top}PD_{l,1})^{-1} ( B_1^{\top}P + \sum_{l=1}^r D_{l,1}^{\top}PC_{lL} + S_{L} ) = 0.
\end{aligned} 
\end{equation}
Throughout this work, $\mathcal{A}$ stands for the set of $L$, where $L \in \mathbb{R}^{m_2 \times n} $ satisfy: 
\begin{itemize}
    \item The system $(A_{L},C_{1L},\dotsb,C_{rL})$ is mean-square stable.
    \item The corresponding ARE \eqref{gtare_lw} has a stabilizing solution $\tilde{P}_{L}$, satisfying the sign conditions
    \[R_{11} + \sum_{l=1}^r D_{l,1}^{\top}\tilde{P}_{L}D_{l,1} \prec 0.\]
\end{itemize}

\subsection{Nested Iterative Scheme}
\label{sec:iterative_design}

We now describe the nested iterative scheme that computes the stabilizing solution of the GTARE \eqref{zslqsdg:gtare}. The scheme requires solving sign-definite AREs as subproblems and generates two sequences:
\begin{itemize}
    \item $\{P^{(k)}\}_{k\ge 0}$: the sequence of successive approximations to the stabilizing solution of the GTARE \eqref{zslqsdg:gtare};
    \item $\{Z^{(k)}\}_{k\ge 0}$: the sequence of semi-definite corrections, each obtained as the stabilizing solution of a sign-definite ARE.
\end{itemize}
The two sequences are linked by $P^{(k+1)} = P^{(k)} + Z^{(k)}$, with $P^{(0)}=0$. At each step, the current iterate $P^{(k)}$ determines the coefficients of a sign-definite ARE whose stabilizing solution provides the correction $Z^{(k)}$; 
% the sign-definiteness of the quadratic term ensures that each sub-problem is tractable via established Riccati solvers.

\paragraph{Initialization ($k=0$).}
Set $P^{(0)} = 0$. Compute $Z^{(0)}$ as the unique stabilizing solution of
\begin{equation}
\label{iterative:initial}
\begin{aligned}
&Z A_{(0)} + A_{(0)}^{\top} Z + \sum_{l=1}^{r}C_{l,(0)}^{\top} Z C_{l,(0)} + M_{(0)}\\
&- \bigl(B_2^{\top} Z + \sum_{l=1}^{r}D_{l,2}^{\top} Z C_{l,(0)}\bigr)^{\top}
\bigl(R_{22} + \sum_{l=1}^{r}D_{l,2}^{\top} Z D_{l,2}\bigr)^{-1}
\bigl(B_2^{\top} Z + \sum_{l=1}^{r}D_{l,2}^{\top} Z C_{l,(0)}\bigr) = 0,
\end{aligned}
\end{equation}
where
\begin{equation}
\label{iterative:matrices_evolve_0}
A_{(0)} = A - \begin{bmatrix} B_1 & B_2 \end{bmatrix} R(0)^{-1}S(0),\,
C_{l,(0)} = C_l - \begin{bmatrix} D_{l,1} & D_{l,2} \end{bmatrix} R(0)^{-1}S(0),\,
M_{(0)} = Q - S^{\top}(0)R(0)^{-1}S(0).
\end{equation}

\paragraph{Iteration step ($k\ge 1$).}
Update
\begin{equation}
\label{iterative:external_circulation}
P^{(k)} = P^{(k-1)} + Z^{(k-1)},
\end{equation}
and compute $Z^{(k)}$ as the unique stabilizing solution of
\begin{equation}
\label{iterative:internal_circulation}
\begin{aligned}
&Z A_{(k)} + A_{(k)}^{\top} Z + \sum_{l=1}^{r}C_{l,(k)}^{\top} Z C_{l,(k)} + M_{(k)}\\
&- \bigl(B_2^{\top}Z + \sum_{l=1}^{r}D_{l,2}^{\top} Z C_{l,(k)}\bigr)^{\top}
\Bigl(R_{22}(P^{(k)}) + \sum_{l=1}^{r}D_{l,2}^{\top} Z D_{l,2}\Bigr)^{-1}
\bigl(B_2^{\top}Z + \sum_{l=1}^{r}D_{l,2}^{\top} Z C_{l,(k)}\bigr)=0,
\end{aligned}
\end{equation}
where the coefficient matrices evolve according to
\begin{equation}
\label{iterative:matrices_evolve_1}
\begin{cases}
A_{(k)} = A_{(k-1)} - \begin{bmatrix} B_1 & B_2 \end{bmatrix} R(P^{(k)})^{-1}\begin{bmatrix} N_{1,(k-1)} \\ N_{2,(k-1)} \end{bmatrix},\\[6pt]
C_{l,(k)} = C_{l,(k-1)} - \begin{bmatrix} D_{l,1} & D_{l,2} \end{bmatrix} R(P^{(k)})^{-1}\begin{bmatrix} N_{1,(k-1)} \\ N_{2,(k-1)} \end{bmatrix},\quad 1\le l\le r,\\[6pt]
M_{(k)} = \bigl(N_{1,(k-1)} - R_{12}(P^{(k)})R_{22}(P^{(k)})^{-1}N_{2,(k-1)}\bigr)^{\top}\\[4pt]
\hspace{2em}\times R^{\sharp}_{22}(P^{(k)})^{-1}
\bigl(N_{1,(k-1)} - R_{12}(P^{(k)})R_{22}(P^{(k)})^{-1}N_{2,(k-1)}\bigr),
\end{cases}
\end{equation}
with auxiliary quantities
\begin{equation}
\label{iterative:matrices_evolve_2}
\begin{aligned}
&R^{\sharp}_{22}(P) = R_{11}(P) - R_{12}(P)R_{22}(P)^{-1}R_{21}(P),\\[4pt]
&\begin{bmatrix} N_{1,(k-1)} \\ N_{2,(k-1)} \end{bmatrix}
=\begin{bmatrix}
B_1^{\top} Z^{(k-1)} + \sum_{l=1}^{r}D_{l,1}^{\top} Z^{(k-1)} C_{l,(k-1)}\\
B_2^{\top} Z^{(k-1)} + \sum_{l=1}^{r}D_{l,2}^{\top} Z^{(k-1)} C_{l,(k-1)}
\end{bmatrix}.
\end{aligned}
\end{equation}

The subsequent analysis establishes two fundamental properties of this scheme:
\begin{itemize}
    \item \textit{Well-posedness.} For every $k\ge 0$, the AREs \eqref{iterative:initial} and \eqref{iterative:internal_circulation} admit unique stabilizing solutions $Z^{(k)}\in \overline{\mathbb{S}}_{+}^n$. Moreover, for any fixed reference strategy $L\in\mathcal{A}$ with associated stabilizing solution $\tilde{P}_L$ of ARE \eqref{gtare_lw}, the uniform bound
    \begin{equation*}
    P^{(k)} = \sum_{s=0}^{k-1} Z^{(s)} \preceq \tilde{P}_L \qquad\text{holds for all } k\ge 0.
    \end{equation*}
    \item \textit{Convergence.} The outer sequence converges: $\lim_{k\to\infty} P^{(k)} = \tilde{P}$, where $\tilde{P}$ is the stabilizing solution of the original GTARE \eqref{zslqsdg:gtare}.
\end{itemize}

\section{Main Results}
\label{sec:main_results}

This section develops the complete theoretical foundation for the nested iterative scheme introduced in Section~\ref{sec:iterative_design}. The analysis is organized in three layers, each building upon the preceding one:

\begin{enumerate}
    \item \textbf{Operator structure} (Section~\ref{subsec:function_g}): We derive algebraic identities for the game operator $\mathcal{G}$ that encode how the GTARE decomposes under perturbations of the solution candidate. These identities---particularly the increment formula of Proposition~\ref{function_g:proposition_3}---are the technical engine driving all subsequent arguments.

    \item \textbf{Subproblem well-posedness} (Section~\ref{subsec:existence}): We establish general conditions under which a sign-definite ARE admits a unique stabilizing solution. This result, encapsulated in Proposition~\ref{main_results:stabilizable_detectable}, guarantees that every inner iteration of our algorithm produces a well-defined correction $Z^{(k)}$.

    \item \textbf{Convergence} (Section~\ref{subsec:convergence}): We prove that the outer sequence $\{P^{(k)}\}_{k\ge 0}$ converges monotonically to the stabilizing solution of the original GTARE \eqref{zslqsdg:gtare}. The proof proceeds by induction, with Lemma~\ref{main_results:lemma} playing a critical role in the inductive argument: it links the system spectral properties across successive iterations. 
\end{enumerate}

\subsection{Structural Properties of the Game Operator}
\label{subsec:function_g}

From the stochastic GTARE \eqref{zslqsdg:gtare} we define a mapping $\mathcal{G}\colon\mathrm{Dom}\,\mathcal{G}\to\mathbb{S}^n$ by
\begin{equation}
\label{main_results:function_g}
    \mathcal{G}(P) = P A + A^{\top} P + \sum_{l=1}^{r}C_{l}^{\top} P C_{l} + Q - S(P)^{\top} R(P)^{-1} S(P).
\end{equation}
The operator $\mathcal{G}$ is not defined on all of $\mathbb{S}^n$; its natural domain is the subset
\[
    \mathrm{Dom}\,\mathcal{G}=\bigl\{P\in\mathbb{S}^n \;\big|\; R_{22}(P) \succ 0 \;\text{ and }\; R_{11}(P) \prec 0\bigr\}.
\]
This operator forms the foundation of our nested iterative scheme. We now derive several of its key properties, which collectively guarantee that the algorithm can be executed iteratively and converges to the desired stabilizing solution.

We first introduce two families of operators that will be used in what follows:
\begin{equation}
\label{main_results:operators_k&n}
\begin{aligned}
\begin{bmatrix} K_1(P) \\ K_2(P) \end{bmatrix} = -R(P)^{-1} \begin{bmatrix} S_1(P) \\ S_2(P) \end{bmatrix},   
\begin{bmatrix}
N_1(P,Z) \\
N_2(P,Z)
\end{bmatrix}=
\begin{bmatrix}
B_1^{\top} Z + \sum_{l=1}^{r} D_{l,1}^{\top} Z \bigl(C_l + D_{l,1}K_1(P) + D_{l,2}K_2(P)\bigr) \\
B_2^{\top} Z + \sum_{l=1}^{r} D_{l,2}^{\top} Z \bigl(C_l + D_{l,1}K_1(P) + D_{l,2}K_2(P)\bigr)
\end{bmatrix},
\end{aligned}
\end{equation}
where $B_i, D_{l,i}$ ($i = 1,2;\, 1 \leq l \leq r$) are defined in \eqref{zslqsdg:sde}, and $R(P)$, $S_1(P)$, $S_2(P)$ are defined in \eqref{zslqsdg:gtare_auxiliarymatrix}.

\begin{proposition}
\label{function_g:proposition_1}
Let $P$ and $Z$ be such that $P, P+Z \in \mathrm{Dom}\,\mathcal{G}$. Then the feedback gains satisfy the following relation:
\[
\begin{bmatrix} K_1(P+Z) \\ K_2(P+Z) \end{bmatrix} 
= \begin{bmatrix} K_1(P) \\ K_2(P) \end{bmatrix}-R(P+Z)^{-1}
\begin{bmatrix}
N_1(P,Z) \\
N_2(P,Z)
\end{bmatrix}.
\]
\end{proposition}

\begin{proof}
Define the increment $\Delta = \begin{bmatrix} K_1(P+Z) \\ K_2(P+Z) \end{bmatrix} - \begin{bmatrix} K_1(P) \\ K_2(P) \end{bmatrix}$. 
By leveraging the definitions of $K_1(\cdot), K_2(\cdot)$ in \eqref{main_results:operators_k&n}, and multiplying both sides of the equation by $R(P+Z)$, we obtain:  
\[
R(P+Z)\Delta = R(P+Z) R(P)^{-1} \begin{bmatrix} S_1(P) \\ S_2(P) \end{bmatrix} - \begin{bmatrix} S_1(P+Z) \\ S_2(P+Z) \end{bmatrix}.
\]
Note that,
\[
R(P+Z) = R(P) + \Delta R,\Delta R = \begin{bmatrix}
\sum_{l=1}^{r} D_{l,1}^{\top} Z D_{l,1} & \sum_{l=1}^{r} D_{l,1}^{\top} Z D_{l,2} \\
\sum_{l=1}^{r} D_{l,2}^{\top} Z D_{l,1} & \sum_{l=1}^{r} D_{l,2}^{\top} Z D_{l,2}
\end{bmatrix}.
\]

Using the identity
$
R(P)^{-1} \begin{bmatrix} S_1(P) \\ S_2(P) \end{bmatrix} = -\begin{bmatrix} K_1(P) \\ K_2(P) \end{bmatrix},
$
we have
\[
R(P+Z)\Delta = \begin{bmatrix} S_1(P) \\ S_2(P) \end{bmatrix} + \Delta R \cdot R(P)^{-1} \begin{bmatrix} S_1(P) \\ S_2(P) \end{bmatrix} - \begin{bmatrix} S_1(P+Z) \\ S_2(P+Z) \end{bmatrix}=-\begin{bmatrix}
N_1(P,Z) \\
N_2(P,Z)
\end{bmatrix}.
\]

Multiplying both sides by $R(P+Z)^{-1}$ yields the relation.
\end{proof}

\begin{proposition}[\cite{Dragan2013book}]
\label{function_g:proposition_2}
For any $P \in \mathrm{Dom}\,\mathcal{G}$ and for all $\varTheta_1 \in \mathbb{R}^{m_1 \times n}, \varTheta_2 \in \mathbb{R}^{m_2 \times n}$, we have
\[
\begin{aligned}
&\mathcal{G}(P) = P (A + B_1 \varTheta_1 + B_2 \varTheta_2) + (A + B_1 \varTheta_1 + B_2 \varTheta_2)^{\top} P + Q\\
& + \sum_{l=1}^{r} (C_{l} + D_{l,1} \varTheta_1 + D_{l,2} \varTheta_2)^{\top} P (C_{l} + D_{l,1} \varTheta_1 + D_{l,2} \varTheta_2)  + \begin{bmatrix} \varTheta_1 \\ \varTheta_2 \end{bmatrix}^{\top} R(0) \begin{bmatrix} \varTheta_1 \\ \varTheta_2\end{bmatrix} \\
&+ \begin{bmatrix} \varTheta_1 \\ \varTheta_2 \end{bmatrix}^{\top} \begin{bmatrix} S_1 \\ S_2 \end{bmatrix} + \begin{bmatrix} S_1 \\ S_2 \end{bmatrix}^{\top} \begin{bmatrix} \varTheta_1 \\ \varTheta_2 \end{bmatrix} - \begin{bmatrix} K_1(P) - \varTheta_1 \\ K_2(P) - \varTheta_2 \end{bmatrix}^{\top} R(P) \begin{bmatrix} K_1(P) - \varTheta_1 \\ K_2(P) - \varTheta_2\end{bmatrix}.
\end{aligned}
\]
\end{proposition}

\begin{proposition}
\label{function_g:proposition_3}
Let $P$, $Z$ satisfy $P, P+Z \in \mathrm{Dom}\,\mathcal{G}$.
Then the following identity holds:
\begin{align*}
&\mathcal{G}(P + Z) = \mathcal{G}(P)+Z(A + B_1K_1(P)+B_2K_2(P))+(A + B_1K_1(P)+B_2K_2(P))^{\top}Z\\
&+\sum_{l=1}^{r}(C_{l} + D_{l,1}K_1(P)+D_{l,2}K_2(P))^{\top}Z(C_{l} + D_{l,1}K_1(P)+D_{l,2}K_2(P)) \\
&- N(P,Z)^{\top} R(P+Z)^{-1} N(P,Z),\quad \text{where} \,\, N(P,Z)=\begin{bmatrix}
N_1(P,Z) \\
N_2(P,Z)
\end{bmatrix}.
\end{align*}
\end{proposition}

\begin{proof}
From the definition of $\mathcal{G}$,
$
\mathcal{G}(P+Z) = (P+Z)A + A^{\top}(P+Z) + \sum_{l=1}^{r} C_l^{\top} (P+Z) C_l + Q - S(P+Z)^{\top} R(P+Z)^{-1} S(P+Z).
$

Applying Proposition \ref{function_g:proposition_2} with $(\varTheta_1, \varTheta_2) = (K_1(P), K_2(P))$ and rearranging, we obtain $\mathcal{G}(P+Z) =$
\[
\begin{aligned}
&P A + A^{\top} P + \sum_{l=1}^{r}C_{l}^{\top} P C_{l}+ Q +K(P)^{\top}R(P)K(P) + S(P)^{\top}K(P)+ K(P)^{\top}S(P)\\
&+Z (A + B_1K_1(P)+B_2K_2(P)) + (A + B_1K_1(P)+B_2K_2(P))^{\top} Z \\
&+ \sum_{l=1}^{r}(C_{l} + D_{l,1}K_1(P)+D_{l,2}K_2(P))^{\top}Z(C_{l} + D_{l,1}K_1(P)+D_{l,2}K_2(P))\\
&- (K(P+Z)-K(P))^{\top} R(P+Z) (K(P+Z)-K(P)),\,\, \text{where} \,\, K(P)=\begin{bmatrix}K_1(P)\\K_2(P)\end{bmatrix}.
\end{aligned}
\]

From Proposition \ref{function_g:proposition_1}, we derive:
\[
(K(P+Z) - K(P))^{\top} R(P+Z) (K(P+Z) - K(P)) = N(P,Z)^{\top} R(P+Z)^{-1} N(P,Z).
\]
Substitute the above term into the rearranged expression of $\mathcal{G}(P+Z)$, and the target identity is thus verified to hold.
\end{proof}

\subsection{Existence and Uniqueness of Stabilizing Solutions for Sign-Definite AREs}
\label{subsec:existence}

The GTARE \eqref{zslqsdg:gtare} involves two distinct types of coupling:
the bias terms appearing in $S(P)$, and the cross-coupling matrix $R_{12}(P)$ appearing quadratically in $R(P)$.
Crucially, these two couplings admit fundamentally different treatments:

\begin{enumerate}
    \item \textit{Bias elimination.}
    For any Riccati equation that arises from a linear-quadratic type cost functional containing state-control cross-weights (such as $S_1$ and $S_2$ here), there exists a \emph{feedback translation}---a constant shift of the control variables---that absorbs these bias terms into the shifted feedback gain.

    \item \textit{Player fixation.}
    Once the bias terms has been eliminated via the feedback translation, imposing a linear relation $v_2(t)=Lx(t)$ on Player~2's shifted control removes the cross-coupling $R_{12}(P)$ from the quadratic term of the GTARE \eqref{zslqsdg:gtare}.
    The reduced problem becomes a single-player stochastic linear-quadratic control problem whose ARE takes the form:
\end{enumerate}
\begin{equation}
\label{pr:are}
\begin{aligned}
&PA + A^{\top}P + \sum_{l=1}^{r} C_{l}^{\top}PC_{l} + Q - \bigl( B^{\top}P + \sum_{l=1}^{r} D_{l}^{\top}PC_{l}  \bigr)^{\top}
\bigl( R + \sum_{l=1}^{r} D_{l}^{\top}PD_{l} \bigr)^{-1}
\bigl( B^{\top}P + \sum_{l=1}^{r} D_{l}^{\top}PC_{l}  \bigr) = 0,
\end{aligned}
\end{equation}
where $B, D_l \in \mathbb{R}^{n \times m}$ for $1 \leq l \leq r$, and $R \in \mathbb{R}^{m \times m}$. The remaining parameters follow the previous definitions.

The decisive feature of ARE \eqref{pr:are} is that its quadratic term involves only the matrix
$R + \sum_{l=1}^{r} D_{l}^{\top}PD_{l}$.
When $R \succ 0$, this ARE falls within the classical theory of sign-definite Riccati equations.
However, an additional subtlety arises: the state-weight matrix $Q$ (e.g., $M_{(k)}$ appearing in \eqref{iterative:internal_circulation}) is typically only positive semidefinite, not strictly positive definite.
Standard sufficient conditions for the existence of stabilizing solutions to sign-definite AREs (e.g., those requiring $Q\succ 0$) are therefore not directly applicable here.
We instead establish general conditions based on mean-square stabilizability and stochastic detectability---assumptions that accommodate the semidefinite case and are natural for the problems at hand---guaranteeing that such an ARE admits a unique stabilizing solution.

We first recall two auxiliary concepts from \cite{Dragan2013book} that will be used in the subsequent analysis.

\begin{definition}[\cite{Dragan2013book}]
The system associated with the ARE \eqref{pr:are} denoted as $[ A, C_1, \dotsb, C_r ; \, B, D_1, \dotsb, D_r ]$ is said to be mean-square stabilizable if there exists a constant matrix $ \varTheta \in\mathbb{R}^{n \times m}$ such that the system $\left(A+B\varTheta, C_{1}+D_1\varTheta, \dotsb, C_{r}+D_r\varTheta\right)$ is mean-square stable.
\end{definition}

\begin{definition}[\cite{Dragan2013book}]
Given the following stochastic observation system:
\begin{equation*}
    \begin{cases}
         dX(t) = A X(t)dt + \sum_{l=1}^{r} C_l X(t)dw_l(t) \\
         dY(t) = E_0 X(t)dt + \sum_{l=1}^{r} E_l X(t)dw_l(t)
    \end{cases},
\end{equation*}
where $E_l \in R^{ q\times n}(0 \leq l \leq r)$, we denote this system by $\left[ E_0, E_1, \dotsb, E_r; A, C_1, \dotsb, C_r \right]$. It is said to be stochastically detectable if there exists a constant matrix $\varTheta \in \mathbb{R}^{n \times q}$ such that the system
$( A + \varTheta E_0,\ C_1 + \varTheta E_1, \dotsb, C_r + \varTheta E_r )$
is mean-square stable.
\end{definition}

\begin{lemma}
\label{main_results:stochastically_detectable_lemma}
If the system $\left[ E_0, E_1, \dotsb, E_r; A, C_1, \dotsb, C_r \right]$ is stochastically detectable, then the following statements are equivalent:
\begin{itemize}
    \item[(a).] The system $(A, C_1, \dotsb, C_r)$ is mean-square stable. Let $\mathcal{L}^*$ denote the linear operator associated with this system, defined by
    $\mathcal{L}^*(P) = P A + A^{\top} P + \sum_{l=1}^{r} C_l^{\top} P C_l; \, \forall P \in \mathbb{S}^n$.
    Then all eigenvalues of $\mathcal{L}^*$ lie in the left half-plane, i.e., $\operatorname{Spec}\mathcal{L}^* \subset \mathbb{C}^{-}$.
    \item[(b).] There exists a matrix $P \in \overline{\mathbb{S}}_{+}^n$ satisfying the matrix equation $P A + A^{\top} P + \sum_{l=1}^{r} C_l^{\top} P C_l + \sum_{l=0}^{r} E_l^{\top} E_l = 0$.
\end{itemize}
\end{lemma}

\begin{proof}
(a) $\Rightarrow$ (b):
Since $\sum_{l=0}^{r} E_l^{\top} E_l \succeq 0$ and the system $(A, C_1, \dotsb, C_r)$ is mean-square stable, it follows from Theorem 3.2.2. and Theorem 2.7.5 in \cite{Dragan2013book} that the equation
$P A + A^{\top} P + \sum_{l=1}^{r} C_l^{\top} P C_l + \sum_{l=0}^{r} E_l^{\top} E_l = 0$
has a solution $P \in \overline{\mathbb{S}}_{+}^n$.

(b) $\Rightarrow$ (a):
This conclusion follows directly from Remark 4.1.5 in \cite{Dragan2013book}. For a detailed argument, see Theorem 4.1.7 and Remark 4.1.5 in the same reference. Furthermore, if the system $(A, C_1, \dotsb, C_r)$ is \textit{mean-square stable}, then $\operatorname{Spec}\mathcal{L}^* \subset \mathbb{C}^{-}$. This is a consequence of Theorem 3.2.2. and Theorem 2.7.7 in \cite{Dragan2013book}.
\end{proof}

Now we define the operator $\boldsymbol{\Lambda}^{[A, C_1, \dotsb, C_r; B, D_1, \dotsb, D_r; Q, R]}: \mathbb{S}^n \rightarrow \mathbb{S}^{n+m}$ associated with the ARE \eqref{pr:are} as 
\begin{equation*}
\boldsymbol{\Lambda}^{[A, C_1, \dotsb, C_r; B, D_1, \dotsb, D_r; Q, R]}(P) =\begin{bmatrix}
PA + A^{\top}P + \sum_{l=1}^{r} C_{l}^{\top}PC_{l} + Q & PB + \sum_{l=1}^{r} C_{l}^{\top}PD_{l} \\
B^{\top}P + \sum_{l=1}^{r} D_{l}^{\top}PC_{l}  & R + \sum_{l=1}^{r} D_{l}^{\top}PD_{l}
\end{bmatrix}.
\end{equation*}

Define the set $\boldsymbol{\Gamma}^{[A, C_1, \dotsb, C_r; B, D_1, \dotsb, D_r; Q, R]}$ related to ARE \eqref{pr:are} as
\begin{equation*}
\left\{ P \in \mathbb{S}^n \,\middle|\, \boldsymbol{\Lambda}^{[A, C_1, \dotsb, C_r; B, D_1, \dotsb, D_r; Q, R]}(P) \succeq 0,\ R + \sum_{l=1}^{r} D_{l}^{\top}PD_{l} \succ 0 \right\}.
\end{equation*}

\begin{proposition}
\label{main_results:stabilizable_detectable}
Suppose the parameters of the ARE \eqref{pr:are} satisfy the following:
\begin{itemize}
    \item $0 \in \boldsymbol{\Gamma}^{[A, C_1, \dotsb, C_r; B, D_1, \dotsb, D_r; Q, R]}$.
    \item The system $[A ,C_1 , \dotsb, C_r ;\ B, D_1, \dotsb, D_r]$ is mean-square stabilizable.
    \item There exists a matrix set $\{E_0, E_1, \dotsb, E_r\}$ satisfying $\sum_{l=0}^{r} E_l^{\top}E_l = Q $ such that the system \[[E_0, E_1, \dotsb, E_r; A, C_1, \dotsb, C_r]\] is stochastically detectable.
\end{itemize}
Then the ARE \eqref{pr:are} admits a unique stabilizing solution $P \in \overline{\mathbb{S}}_{+}^n$ such that $R + \sum_{l=1}^{r} D_{l}^{\top}PD_{l} \succ 0$.
\end{proposition}

\begin{proof}
By Theorem 4.7 in \cite{Dragan2004}, along with the condition $0 \in \boldsymbol{\Gamma}^{[A, C_1, \dotsb, C_r; B, D_1, \dotsb, D_r; Q, R]}$ and the mean-square stabilizability of system $\left[ A , C_1 , \dotsb, C_r ; B, D_1, \dotsb, D_r \right]$,   
the ARE \eqref{pr:are} has a maximal solution $P_{\text{max}} \in \overline{\mathbb{S}}_{+}^n$ such that $R + \sum_{l=1}^{r} D_{l}^{\top}P_{\text{max}}D_{l} \succ 0$. 

Using Proposition \ref{function_g:proposition_2}, the ARE for $P_{\text{max}}$ can be rewritten as
\[
\begin{aligned}
&P_{\text{max}}\left( A  + B T(P_{\text{max}}) \right) + \left( A  + B T(P_{\text{max}}) \right)^{\top} P_{\text{max}} \\
&+ \sum_{l=1}^{r} \left(C_l  + D_l T(P_{\text{max}})\right)^{\top} P_{\text{max}} \left(C_l  + D_l T(P_{\text{max}})\right) + \sum_{l=0}^{r} \hat{E}_l^{\top} \hat{E}_l = 0,
\end{aligned}
\]
where $ T(P) = - \left( R + \sum_{l=1}^{r} D_{l}^{\top}PD_{l} \right)^{-1} \left( B^{\top}P + \sum_{l=1}^{r} D_{l}^{\top}PC_{l} \right)$, and
\[
\hat{E}_l = 
\begin{bmatrix}
\mathbb{O}_{q \times n}^{\times (l)} \\
E_l \\
\mathbb{O}_{q \times n}^{\times (r-l)} \\
\mathbb{O}_{m \times n}^{\times (l)} \\
\sqrt{\frac{1}{r+1}R}\,T(P_{\text{max}}) \\
\mathbb{O}_{m \times n}^{\times (r-l)}
\end{bmatrix} \in \mathbb{R}^{[(q+m)(r+1)] \times n},\,0 \leq l \leq r .
\]

Since the system $\left[E_0, E_1, \dotsb, E_r; A , C_1 , \dotsb, C_r \right]$ is stochastically detectable, there exists $\varTheta \in \mathbb{R}^{n \times q}$ such that the system
$
\left(A  + \varTheta E_0, C_1  + \varTheta E_1, \dotsb, C_r + \varTheta E_r \right)
$ 
is mean-square stable.
Let $\hat{\varTheta}\in \mathbb{R}^{n \times [(q+m)(r+1)]}$ and $\hat{\varTheta} =$  
\[
\begin{bmatrix}
\varTheta^{\times (r+1)} & -B\sqrt{(r+1)R^{-1}} & -D_1\sqrt{(r+1)R^{-1}} & \dotsb & -D_r\sqrt{(r+1)R^{-1}}
\end{bmatrix},
\]  
then the system  
\[
\left(A  + B T(P_{\text{max}}) + \hat{\varTheta} \hat{E}_0, C_1 + D_1T(P_{\text{max}}) + \hat{\varTheta} \hat{E}_1, \dotsb, C_r + D_rT(P_{\text{max}}) + \hat{\varTheta} \hat{E}_r \right)
\]  
is mean-square stable. This means the system
\[
\left[ \hat{E}_0, \hat{E}_1, \dotsb, \hat{E}_r;A  + B T(P_{\text{max}}), C_1 + D_1T(P_{\text{max}}), \dotsb, C_r + D_rT(P_{\text{max}}) \right]
\]
is stochastically detectable.

By Lemma \ref{main_results:stochastically_detectable_lemma}, we get the system  
\[
(A + BT(P_{\text{max}}),\ C_1 + D_1T(P_{\text{max}}), \dotsb, C_r + D_rT(P_{\text{max}}))
\]  
is mean-square stable. Hence, $P_{\text{max}}$ is a stabilizing solution of ARE \eqref{pr:are}. By Theorem 5.6.5 in \cite{Dragan2013book}, the stabilizing solution is unique, so ARE \eqref{pr:are} has a unique stabilizing solution $P\in \overline{\mathbb{S}}_{+}^n$ such that $R + \sum_{l=1}^{r} D_{l}^{\top}PD_{l} \succ 0$.
\end{proof}

\subsection{Convergence Analysis of the Iterative Sequence}
\label{subsec:convergence}

The convergence proof proceeds by induction, relying on previously derived structural properties. Proposition~\ref{main_results:stabilizable_detectable} ensures a unique stabilizing solution for each sign-definite ARE arising in the inner iteration, rendering the sequences $\{Z^{(k)}\}_{k\ge 0}$ and $\{P^{(k)}\}_{k\ge 0}$ well-defined. As the inductive cornerstone, Lemma~\ref{main_results:lemma} builds a two-way link between the spectral stability of the closed-loop system at the current iterate and the partial order among stabilizing solutions of successive iterates. The monotonicity and boundedness of the sequences then guarantee convergence to the stabilizing solution of GTARE \eqref{zslqsdg:gtare}.

For each triple $(P,Z,L)$ satisfying $P\in\mathrm{Dom}\,\mathcal{G}$, $P + Z\in\mathrm{Dom}\,\mathcal{G}$ and $L \in \mathcal{A}$, the linear operators $\mathcal{L}^{*}_{P}$ and $\mathcal{L}^{*}_{P+Z}$ are defined as follows:
\begin{equation}
\label{linear_operators_pl}
\begin{aligned}
    &\mathcal{L}^{*}_{P}(Y) = Y ( A_{(0)}+ B_1\hat{K}_1(P)+B_2L)+( A_{(0)}+ B_1\hat{K}_1(P)+B_2L )^{\top}Y +\\
    &\sum_{l=1}^r (C_{l,(0)}+D_{l,1}\hat{K}_1(P)+D_{l,2}L)^{\top} Y (C_{l,(0)}+D_{l,1}\hat{K}_1(P)+D_{l,2}L), \, \forall \,Y \in \mathbb{S}^n;
\end{aligned}
\end{equation}
\begin{equation}
\label{linear_operators_pzl}
\begin{aligned}
    &\mathcal{L}^{*}_{P+Z}(Y) = Y ( A_{(0)}+ B_1\hat{K}_1(P+Z)+B_2L)+( A_{(0)}+ B_1\hat{K}_1(P+Z)+B_2L)^{\top} Y +\\
    & \sum_{l=1}^r (C_{l,(0)}+ D_{l,1}\hat{K}_1(P+Z)+D_{l,2}L)^{\top} Y (C_{l,(0)}+D_{l,1}\hat{K}_1(P+Z)+D_{l,2}L), \, \forall \, Y \in \mathbb{S}^n,
\end{aligned}
\end{equation}
where $A_{(0)},C_{l,(0)}(1 \leq l \leq r)$ are defined in \eqref{iterative:matrices_evolve_0} and
\begin{equation}
\label{linearoperators_auxiliarymatrix}
\begin{bmatrix} \hat{K}_1(P) \\ \hat{K}_2(P) \end{bmatrix}
= \begin{bmatrix} K_1(P) \\ K_2(P) \end{bmatrix}-\begin{bmatrix} K_1(0) \\ K_2(0) \end{bmatrix},\quad \begin{bmatrix} K_1(P) \\ K_2(P) \end{bmatrix} = -R(P)^{-1} \begin{bmatrix} S_1(P) \\ S_2(P) \end{bmatrix}.
\end{equation}

\begin{lemma}
\label{main_results:lemma}
Assume $0 \in \mathrm{Dom}\,\mathcal{G}$ and there exists an $L \in \mathcal{A}$. Let $P, Z\in \mathbb{S}^n$ satisfy the following conditions:
\begin{itemize}
    \item $P, P+Z \in \mathrm{Dom}\,\mathcal{G}$;
    \item The matrices $P$ and $Z$ satisfy the algebraic relation
    \begin{equation}
    \label{main_results:g_interrelated_equation}
    \begin{aligned}
        0=&\mathcal{G}(P)+Z\left(A_{(0)} + B_1\hat{K}_1(P) + B_2\hat{K}_2(P) \right)+\left(A_{(0)} +  B_1\hat{K}_1(P) + B_2\hat{K}_2(P)\right)^{\top}Z\\
        &+\sum_{l=1}^{r}\left(C_{l,(0)} +  D_{l,1}\hat{K}_1(P) + D_{l,2}\hat{K}_2(P)\right)^{\top}Z\left(C_{l,(0)} +  D_{l,1}\hat{K}_1(P) + D_{l,2}\hat{K}_2(P)\right)\\
        &-\left[B_2^{\top} Z + \sum_{l=1}^{r}D_{l,2}^{\top} Z \left(C_{l,(0)} +  D_{l,1}\hat{K}_1(P) + D_{l,2}\hat{K}_2(P)\right)\right]^{\top} R_{22}(P+Z)^{-1}\\
        & \quad \times \left[B_2^{\top} Z + \sum_{l=1}^{r}D_{l,2}^{\top} Z \left(C_{l,(0)} +  D_{l,1}\hat{K}_1(P) + D_{l,2}\hat{K}_2(P)\right)\right].
    \end{aligned}
    \end{equation}
\end{itemize}
Let $\mathcal{L}^*_{P}$ and $\mathcal{L}^*_{P+Z}$ denote the linear operators associated with the triple $(P,Z,L)$ defined by \eqref{linear_operators_pl} and \eqref{linear_operators_pzl}, respectively. Suppose that $\tilde{P}_{L}$ is the stabilizing solution of the ARE \eqref{gtare_lw} corresponding to $L$. Then the following conclusions hold:
\begin{enumerate}
    \item[(\romannumeral1)] If $\operatorname{Spec} \,\mathcal{L}^*_{P}\subset\mathbb{C}^{-}$, then $\tilde{P}_{L}\succeq P + Z$;
    \item[(\romannumeral2)] If $\tilde{P}_{L}\succeq P + Z$, then $\operatorname{Spec}\, \mathcal{L}^*_{P+Z}\subset\mathbb{C}^{-}$.
\end{enumerate}
\end{lemma}

\begin{proof}
We first prove conclusion (\romannumeral1). Combining Proposition \ref{function_g:proposition_2} and Proposition \ref{function_g:proposition_3}, we derive that
\begin{align*}
    \mathcal{G}(P+Z) =&\mathcal{L}^*_{P}(P+Z) + Q_{L}+\hat{K}^{\top}_1(P)R_{12}L +L^{\top}R_{21}\hat{K}_1(P)+ \hat{K}_1^{\top}(P)R_{11}\hat{K}_1(P)\\
    &-\left(J(P)-\hat{K}(P+Z)\right)^{\top}R(P+Z)\left(J(P)-\hat{K}(P+Z)\right),
\end{align*}
where $J(P)=\begin{bmatrix}\hat{K}_1(P) \\ L\end{bmatrix}$, the linear operator $\mathcal{L}^{*}_{P}$ is defined in \eqref{linear_operators_pl}, and the mappings $\hat{K}_1(P)$ and $\hat{K}(P)$ are given in \eqref{linearoperators_auxiliarymatrix}.

By virtue of \eqref{main_results:g_interrelated_equation} and Proposition \ref{function_g:proposition_3}, it can be verified that $P + Z$ is a solution to the following ARE:
\begin{equation}
\label{lemma_proof_1}
\begin{aligned}
    &\mathcal{L}^*_{P}(P+Z) + Q_{L}+\hat{K}^{\top}_1(P)R_{12}L +L^{\top}R_{21}\hat{K}_1(P)+ \hat{K}_1^{\top}(P)R_{11}\hat{K}_1(P)\\
    &-\left(J(P)-\hat{K}(P+Z)\right)^{\top}R(P+Z)\left(J(P)-\hat{K}(P+Z)\right)\\
    &+\left( \hat{N}_1(P,Z)-R_{12}(P+Z)R_{22}(P+Z)^{-1}\hat{N}_2(P,Z) \right)^{\top}\\
    & \quad \times R_{22}^{\sharp}(P + Z)^{-1}\left( \hat{N}_1(P,Z)-R_{12}(P+Z)R_{22}(P+Z)^{-1}\hat{N}_2(P,Z) \right)=0,
\end{aligned}
\end{equation}
where
\begin{align*}
    \hat{N}_1(P,Z)&=B_1^{\top} Z + \sum_{l=1}^{r}D_{l,1}^{\top} Z \left(C_{l,(0)} + D_{l,1}\hat{K}_1(P)+D_{l,2}\hat{K}_2(P)\right),\\
    \hat{N}_2(P,Z)&=B_2^{\top} Z + \sum_{l=1}^{r}D_{l,2}^{\top} Z \left(C_{l,(0)} + D_{l,1}\hat{K}_1(P)+D_{l,2}\hat{K}_2(P)\right),
\end{align*}
and the matrix $R_{22}^{\sharp}(P + Z)$ is defined in \eqref{iterative:matrices_evolve_2}. Furthermore, we have the decomposition
\begin{equation}
\label{lemma_proof_2}
\begin{aligned}
&\left(J(P)-\hat{K}(P+Z)\right)^{\top}R(P+Z)\left(J(P)-\hat{K}(P+Z)\right)=\\
&H^{\top}_1H_1+\left(\hat{K}_1(P)- \hat{K}_1(P+Z) \right)^\top R_{22}^{\sharp}(P+Z)\left(\hat{K}_1(P)- \hat{K}_1(P+Z) \right)=\\
& H^{\top}_1H_1+\left( \hat{N}_1(P,Z)-R_{12}(P+Z)R_{22}(P+Z)^{-1}\hat{N}_2(P,Z) \right)^{\top}R_{22}^{\sharp}(P + Z)^{-1}\\
& \quad \times \left( \hat{N}_1(P,Z)-R_{12}(P+Z)R_{22}(P+Z)^{-1}\hat{N}_2(P,Z) \right),
\end{aligned}
\end{equation}
where $H_1=\begin{bmatrix}R_{22}(P + Z)^{-\frac{1}{2}}R_{21}(P + Z) & R_{22}(P + Z)^{\frac{1}{2}}\end{bmatrix} \left( J(P)-\hat{K}(P+Z)\right)$.

According to Proposition \ref{function_g:proposition_2}, the stabilizing solution $\tilde{P}_{L}$ satisfies the equation
\begin{equation}
\label{lemma_proof_3}
\begin{aligned}
&\mathcal{L}^{*}_{P}(\tilde{P}_{L})+ Q_{L}+\hat{K}^{\top}_1(P)R_{12}L +L^{\top}R_{21}\hat{K}_1(P)+ \hat{K}_1^{\top}(P)R_{11}\hat{K}_1(P)\\
&- \left( \hat{K}_1(P) - K_{L}(\tilde{P}_{L})\right)^\top (R_{11} + \sum_{l=1}^r D_{l,1}^{\top}\tilde{P}_{L}D_{l,1}) \left( \hat{K}_1(P) - K_{L}(\tilde{P}_{L})\right)=0,
\end{aligned} 
\end{equation}
with the gain matrix given by
\[
K_{L}(\tilde{P}_{L})=-(R_{11} + \sum_{l=1}^r D_{l,1}^{\top}\tilde{P}_{L}D_{l,1})^{-1}(B^{\top}_1\tilde{P}_{L} + \sum_{l=1}^r D^{\top}_{l,1}\tilde{P}_{L}C_{l} + R_{12}L).
\]
Subtracting \eqref{lemma_proof_1} from \eqref{lemma_proof_3} and incorporating the decomposition \eqref{lemma_proof_2}, we obtain
\begin{align*}
    &\mathcal{L}^{*}_{P}(\tilde{P}_{L}-P- Z) +H^{\top}H=0,\\
    &H^{\top}H\triangleq H_1^{\top}H_1- \left( \hat{K}_1(P) - K_{L}(\tilde{P}_{L})\right)^\top (R_{11} + \sum_{l=1}^r D_{l,1}^{\top}\tilde{P}_{L}D_{l,1}) \left( \hat{K}_1(P) - K_{L}(\tilde{P}_{L})\right).   
\end{align*}
Since $\operatorname{Spec} \,\mathcal{L}^*_{P}\subset\mathbb{C}^{-}$ and $H^{\top}H\succeq0$, analogous to the proof of implication (a) $\Rightarrow$ (b) in Lemma \ref{main_results:stochastically_detectable_lemma}, we conclude that $\tilde{P}_{L} \succeq P + Z$. This completes the proof of conclusion (\romannumeral1).

We now proceed to prove conclusion (\romannumeral2). Similarly, by applying Proposition \ref{function_g:proposition_2} and Proposition \ref{function_g:proposition_3} along with condition \eqref{main_results:g_interrelated_equation}, $P + Z$ satisfies the ARE
\begin{equation}
\label{lemma_proof_4}
\begin{aligned}
    &\mathcal{L}^*_{P+Z}(P+Z) +\hat{K}^{\top}_1(P+Z)R_{12}L +L^{\top}R_{21}\hat{K}_1(P+Z)+\hat{K}_1^{\top}(P+Z)R_{11}\hat{K}_1(P+Z)\\
    &+ Q_{L}-\left(J(P+Z)-\hat{K}(P+Z)\right)^{\top}R(P+Z)\left(J(P+Z)-\hat{K}(P+Z)\right)\\
    &+\left(\hat{N}_1(P,Z)-R_{12}(P+Z)R_{22}(P+Z)^{-1}\hat{N}_2(P,Z) \right)^{\top}R_{22}^{\sharp}(P + Z)^{-1}\\
    &\quad \times\left(\hat{N}_1(P,Z)-R_{12}(P+Z)R_{22}(P+Z)^{-1}\hat{N}_2(P,Z) \right)=0.
\end{aligned}
\end{equation}
In addition, it follows from Proposition \ref{function_g:proposition_2} that $\tilde{P}_{L}$ satisfies the equation
\begin{equation}
\label{lemma_proof_5}
\begin{aligned}
&\mathcal{L}^{*}_{P+Z}(\tilde{P}_{L})+\hat{K}^{\top}_1(P+Z)R_{12}L +L^{\top}R_{21}\hat{K}_1(P+Z)+ \hat{K}_1^{\top}(P+Z)R_{11}\hat{K}_1(P+Z)+\\
&Q_{L}- \left( \hat{K}_1(P+Z) - K_{L}(\tilde{P}_{L})\right)^\top (R_{11} + \sum_{l=1}^r D_{l,1}^{\top}\tilde{P}_{L}D_{l,1}) \left( \hat{K}_1(P+Z) - K_{L}(\tilde{P}_{L})\right)=0.
\end{aligned} 
\end{equation}
Moreover, the quadratic term of \eqref{lemma_proof_4} can be decomposed as
\begin{equation}
\label{lemma_proof_6}
\begin{aligned}
&\left(J(P+Z)-\hat{K}(P+Z)\right)^{\top}R(P+Z)\left(J(P+Z)-\hat{K}(P+Z)\right)=\\
&\left(\hat{K}_1(P+Z)- \hat{K}_1(P+Z) \right)^\top R_{22}^{\sharp}(P+Z)\left(\hat{K}_1(P+Z)- \hat{K}_1(P+Z) \right)+H^{\top}_2H_2\\
&=H^{\top}_2H_2,
\end{aligned}
\end{equation}
where $H_2=\begin{bmatrix}R_{22}(P + Z)^{-\frac{1}{2}}R_{21}(P + Z) & R_{22}(P + Z)^{\frac{1}{2}}\end{bmatrix} \left( J(P+Z)-\hat{K}(P+Z)\right)$.

Subtracting \eqref{lemma_proof_4} from \eqref{lemma_proof_5} and combining decomposition \eqref{lemma_proof_6}, we derive
\begin{align*}
    & \mathcal{L}^{*}_{P+Z}(\tilde{P}_{L}-P - Z) + H_2^{\top}H_2- \left( \hat{K}_1(P+Z) - K_{L}(\tilde{P}_{L})\right)^\top (R_{11} + \sum_{l=1}^r D_{l,1}^{\top}\tilde{P}_{L}D_{l,1}) \\
    &\times\left(\hat{K}_1(P+Z) - K_{L}(\tilde{P}_{L})\right)- \left(\hat{N}_1(P,Z)-R_{12}(P+Z)R_{22}(P+Z)^{-1}\hat{N}_2(P,Z) \right)^{\top}\\
    &\times R_{22}^{\sharp}(P + Z)^{-1}\left(\hat{N}_1(P,Z)-R_{12}(P+Z)R_{22}(P+Z)^{-1}\hat{N}_2(P,Z) \right) =0.  
\end{align*}
Define the matrix deviation $\Delta\triangleq\tilde{P}_{L}-P - Z$. Then we have $\mathcal{L}^{*}_{P+Z}(\Delta)+\sum_{l=0}^{r}\hat{E}_l^{\top}\hat{E}_l=0$ with $\Delta \in \overline{\mathbb{S}}_{+}^n$, where $\hat{E}_l\in\mathbb{R}^{[m_2 + m_1(r+2)]\times n}$ for $0 \leq l \leq r$, and
\[
    \hat{E}_l=\begin{bmatrix}  \mathbb{O}_{m_1 \times n}^{\times (l)}\\\sqrt{-\frac{1}{r+1}\left(R_{11} + \sum_{l=1}^r D_{l,1}^{\top}\tilde{P}_{L} D_{l,1}\right)} \left(\hat{K}_1(P+Z) - K_{L}(\tilde{P}_{L})\right)\\\mathbb{O}_{m_1 \times n}^{\times (r-l)}\\ \frac{1}{\sqrt{r+1}}H_2 \\ \sqrt{-\frac{1}{r+1}R_{22}^{\sharp}(P+ Z)^{-1}}\left(\hat{N}_1(P,Z)-R_{12}(P+Z)R_{22}(P+Z)^{-1}\hat{N}_2(P,Z)\right)\end{bmatrix}.
\]
Introduce the auxiliary matrix $\hat{\varTheta} \in\mathbb{R}^{n\times [m_2 + m_1(r+2)]}$ such that
\[
\hat{\varTheta}=\begin{bmatrix}  -\sqrt{-(r+1)\left(R_{11} + \sum_{l=1}^r D_{l,1}^{\top}\tilde{P}_{L} D_{l,1}\right)^{-1}}B_{1}^{\top}\\-\sqrt{-(r+1)\left(R_{11} + \sum_{l=1}^r D_{l,1}^{\top}\tilde{P}_{L} D_{l,1}\right)^{-1}}D_{1,1}^{\top}\\ \dotsb \\-\sqrt{-(r+1)\left(R_{11} + \sum_{l=1}^r D_{l,1}^{\top}\tilde{P}_{L} D_{l,1}\right)^{-1}}D_{r,1}^{\top} \\ \mathbb{O}_{n \times m_2}^{\top}\\ \mathbb{O}_{n \times m_1}^{\top} \end{bmatrix}^{\top}.
\]
The system
$ 
(A_{(0)}+ B_1\hat{K}_1(P+Z)+B_2L+\hat{\varTheta}\hat{E}_0,C_{1,(0)}+D_{1,1}\hat{K}_1(P+Z)+D_{1,2}L+\hat{\varTheta}\hat{E}_1,\dotsb,C_{r,(0)}+D_{r,1}\hat{K}_1(P+Z)+D_{r,2}L+\hat{\varTheta}\hat{E}_r) 
$
associated with the system 
$
(A_{L}+B_{1}K_{L}(\tilde{P}_{L}),C_{1L}+D_{1,1}K_{L}(\tilde{P}_{L}),\dotsb,C_{rL}+D_{r,1}K_{L}(\tilde{P}_{L})) 
$ is mean-square stable.
This implies that the system
$ 
[\hat{E}_0,\hat{E}_1,\dotsb,\hat{E}_r;A_{(0)}+ B_1\hat{K}_1(P+Z)+B_2L,C_{1,(0)}+D_{1,1}\hat{K}_1(P+Z)+D_{1,2}L,\dotsb,C_{r,(0)}+D_{r,1}\hat{K}_1(P+Z)+D_{r,2}L]
$  is stochastically detectable. By virtue of Lemma \ref{main_results:stochastically_detectable_lemma}, we conclude that $\operatorname{Spec}\, \mathcal{L}^*_{P+Z} \subset \mathbb{C}^{-}$. The proof of conclusion (\romannumeral2) is therefore completed.
\end{proof}

Before stating the subsequent theorem, we first introduce two classes of linear operators $\mathcal{L}^{*}_{P^{(k)}}$ and $\mathcal{L}^{(k,k+1)}$ with $k = 0,1,2,\dots$, which are associated with the iterative sequences $\{Z^{(k)}\}_{k\geq0}$ and $\{P^{(k)}\}_{k\geq0}$ defined in Section~\ref{sec:iterative_design}.
\begin{equation}
\label{linear_operators_pkl}  
\begin{aligned}
&\mathcal{L}^{*}_{P^{(k)}}(Y) = \, Y \left( A_{(0)}+B_1\hat{K}_1(P^{(k)})+B_2L\right)+\left( A_{(0)}+B_1\hat{K}_1(P^{(k)})+B_2L\right)^{\top}Y \\
&+\sum_{l=1}^r \left(C_{l,(0)}+D_{l,1}\hat{K}_1(P^{(k)})+D_{l,2}L\right)^{\top} Y \left(C_{l,(0)}+D_{l,1}\hat{K}_1(P^{(k)})+D_{l,2}L\right),\, \forall\, Y \in \mathbb{S}^n;
\end{aligned} 
\end{equation}
\begin{equation}
\label{linear_operators_k&+1}
\begin{aligned}
&\mathcal{L}^{(k,k+1)}(Y) = \, Y \left( A_{(k)}+B_2T_{(k+1)}\right) + \left( A_{(k)}+B_2T_{(k+1)}\right)^{\top} Y \\
&+ \sum_{l=1}^r \left(  C_{l,(k)}+D_{l,2}T_{(k+1)}\right)^{\top} Y \left(  C_{l,(k)}+D_{l,2}T_{(k+1)}\right),\, \forall \,Y \in \mathbb{S}^n,
\end{aligned}
\end{equation}
where $T_{(k+1)}=-R_{22}(P^{(k)}+Z^{(k)})^{-1}N_{2,(k)}$. The matrices $A_{(k)}$, $C_{l,(k)}(1 \leq l \leq r)$ and $N_{2,(k)}$ are defined in \eqref{iterative:matrices_evolve_1} and \eqref{iterative:matrices_evolve_2}.

\begin{theorem}
\label{main_results:theorem}
Suppose the following assumptions hold:
\begin{itemize}
    \item $R_{22} \succ 0$ and there exists $L \in \mathcal{A}$;
    \item There exist matrices $E_{l,(0)}$ for $0 \leq l \leq r$ such that $\sum_{l=0}^{r} E_{l,(0)}^\top E_{l,(0)} = Q - S^\top(0) R(0)^{-1}S(0)$, and the system
    \[
    \left[E_{0,(0)}, E_{1,(0)}, \dotsb, E_{r,(0)}; A_{(0)}, C_{1,(0)}, \dotsb, C_{r,(0)}\right]
    \]
    is stochastically detectable.
\end{itemize}
Then the following conclusions hold:

\begin{enumerate}
    \item The sequences $\{Z^{(k)}\}_{k\geq0}$ and $\{P^{(k)}\}_{k\geq0}$ are well-defined via \eqref{iterative:initial}, \eqref{iterative:external_circulation} and \eqref{iterative:internal_circulation}. For $k = 0,1,2,\dots$, the following properties are satisfied:
    \begin{enumerate}
        \item[$a_k$] $\operatorname{Spec}\, \mathcal{L}^{(k,k+1)}\subset\mathbb{C}^{-}$;
        \item[$b_k$] Fix an arbitrary $L \in \mathcal{A}$, and let $\tilde{P}_{L}$ be the stabilizing solution of the corresponding ARE \eqref{gtare_lw}. Then $\tilde{P}_{L} \succeq P^{(k)} + Z^{(k)}$ and $\operatorname{Spec}\, \mathcal{L}^*_{P^{(k+1)}}\subset\mathbb{C}^{-}$;
        \item[$c_k$] $P^{(k)},P^{(k)}+Z^{(k)} \in \mathrm{Dom}\,\mathcal{G}$, and
        \begin{align*}
        \mathcal{G}(P^{(k)}+Z^{(k)})
        =&- \left( N_{1,(k)}-R_{12}(P^{(k)}+Z^{(k)})R_{22}(P^{(k)}+Z^{(k)})^{-1}N_{2,(k)} \right)^{\top} \\
        &\times R_{22}^{\sharp}(P^{(k)}+Z^{(k)})^{-1} \\
        &\times\left( N_{1,(k)}-R_{12}(P^{(k)}+Z^{(k)})R_{22}(P^{(k)}+Z^{(k)})^{-1}N_{2,(k)} \right),   
        \end{align*}
        where $R_{22}^{\sharp}(P^{(k)}+Z^{(k)})$, $N_{1,(k)}$ and $N_{2,(k)}$ are given in \eqref{iterative:matrices_evolve_2}.
    \end{enumerate}
    
    \item $\lim\limits_{k\rightarrow\infty}P^{(k)}=\tilde{P}$, where $\tilde{P}$ denotes the unique stabilizing solution of the stochastic GTARE \eqref{zslqsdg:gtare}.
\end{enumerate}
\end{theorem}

\begin{proof}
\textbf{Step 1. Initial case $k=0$.}
For $k = 0$, we have $P^{(0)} = 0$, and $Z^{(0)}$ is the stabilizing solution of the following algebraic Riccati equation:
\begin{equation}
\label{proof_1}
\begin{aligned}
&ZA_{(0)} + A_{(0)}^{\top} Z + \sum_{l=1}^{r}C_{l,(0)}^{\top} ZC_{l,(0)} + Q - S^\top(0) R(0)^{-1}S(0)\\
&- \bigl(B_2^{\top} Z + \sum_{l=1}^{r}D_{l,2}^{\top} Z C_{l,(0)} \bigr)^{\top}
\bigl(R_{22} + \sum_{l=1}^{r}D_{l,2}^{\top} Z D_{l,2}\bigr)^{-1}
\bigl(B_2^{\top} Z + \sum_{l=1}^{r}D_{l,2}^{\top} Z C_{l,(0)}\bigr) = 0.    
\end{aligned}    
\end{equation}
We first verify the well-posedness of $Z^{(0)}$ by establishing the existence and uniqueness of its stabilizing solution.

Since $\mathcal{A}$ is nonempty, take an arbitrary fixed $L \in \mathcal{A}$ and set $L^{(0)} = L$. One can verify that the system $(A_{(0)}+B_{2}L^{(0)},C_{1,(0)}+D_{1,2}L^{(0)},\dotsb,C_{r,(0)}+D_{r,2}L^{(0)})$ associated with $(A_L,C_{1L},\dotsb,C_{rL})$ is mean-square stable, which implies that the system
\[
[A_{(0)},C_{1,(0)},\dotsb,C_{r,(0)};B_2,D_{1,2},\dotsb,D_{r,2}]
\]
is mean-square stabilizable. In addition, the condition $R_{22} \succ 0$ together with 
\[
\sum_{l=0}^{r}E^{\top}_{l,(0)}E_{l,(0)} = Q - S^\top(0) R(0)^{-1}S(0)\succeq 0
\]
yields
\[
0 \in \boldsymbol{\Gamma}^{[A_{(0)},C_{1,(0)},\dotsb,C_{r,(0)};B_2,D_{1,2},\dotsb,D_{r,2};\sum_{l=0}^{r}E^{\top}_{l,(0)}E_{l,(0)},R_{22}]}.
\]
Combined with the stochastic detectability of the system $[E_{0,(0)},E_{1,(0)},\dotsb,E_{r,(0)}; A_{(0)},C_{1,(0)},\dotsb,C_{r,(0)}]$, it follows from Proposition \ref{main_results:stabilizable_detectable} that ARE \eqref{proof_1} admits a unique stabilizing solution $Z^{(0)}$ satisfying $R_{22}(Z^{(0)}) \succ 0$. Hence $Z^{(0)}$ is well-defined, which further implies $\operatorname{Spec}\, \mathcal{L}^{(0,1)}\subset\mathbb{C}^{-}$.

Let $\tilde{P}_{L}$ be the stabilizing solution of ARE \eqref{gtare_lw} corresponding to $L$. Since the operator $\mathcal{L}^*_{P^{(0)}}$ is associated with the system $(A_{L},C_{1L},\dotsb,C_{rL})$ is mean-square stable, we have the $\operatorname{Spec}\, \mathcal{L}^*_{P^{(0)}}\subset\mathbb{C}^{-}$. Applying conclusions (i) and (ii) of Lemma \ref{main_results:lemma}, we obtain $\tilde{P}_{L} \succeq P^{(0)} + Z^{(0)}$ and $\operatorname{Spec}\, \mathcal{L}^*_{P^{(1)}}\subset\mathbb{C}^{-}$.

From $R_{22}(\tilde{P}_{L}) \succeq R_{22} \succ 0$ and $R_{11}\preceq R_{11}(\tilde{P}_{L}) \prec 0$, we know $P^{(0)},\tilde{P}_{L} \in \mathrm{Dom}\,\mathcal{G}$. Meanwhile, $R_{22}(P^{(0)} + Z^{(0)})\succeq R_{22} \succ 0$ and $R_{11}(P^{(0)} + Z^{(0)}) \prec R_{11}(\tilde{P}_{L}) \prec 0$ ensure $P^{(0)} + Z^{(0)} \in \mathrm{Dom}\,\mathcal{G}$.
Substituting \eqref{proof_1} into the expression of $\mathcal{G}( P^{(0)} + Z^{(0)} )$ and utilizing Proposition \ref{function_g:proposition_3}, we deduce
\[
\mathcal{G}( P^{(0)} + Z^{(0)} ) 
=-\left(N_{1,(0)}-R_{12}(Z^{(0)})R_{22}(Z^{(0)})^{-1}N_{2,(0)}\right)^{\top}R^{\sharp}_{22}(Z^{(0)})^{-1}\left(N_{1,(0)}-R_{12}(Z^{(0)})R_{22}(Z^{(0)})^{-1}N_{2,(0)}\right).
\]
The proof of $a_0$--$c_0$ is thus completed.

\textbf{Step 2. Inductive step.}
Assume that for $k=h-1$, $Z^{(h-1)}$ is well-defined and all properties $a_{h-1}$--$c_{h-1}$ hold. We proceed to prove that $Z^{(h)}$ is well-defined and $a_{h}$--$c_{h}$ are valid.

For $k=h$, $Z^{(h)}$ satisfies the following ARE:
\begin{equation}
\label{proof_2}
\begin{aligned}
&ZA_{(h)} + A_{(h)}^{\top} Z + \sum_{l=1}^{r}C_{l,(h)}^{\top} ZC_{l,(h)} + \sum_{l=0}^{r}E^{\top}_{l,(h)}E_{l,(h)}\\
&- \bigl( B_2^{\top} Z + \sum_{l=1}^{r}D_{l,2}^{\top} Z C_{l,(h)} \bigr)^{\top}
\Bigl( R_{22}(P^{(h)}) + \sum_{l=1}^{r}D_{l,2}^{\top} Z D_{l,2} \Bigr)^{-1}
\bigl( B_2^{\top} Z + \sum_{l=1}^{r}D_{l,2}^{\top} Z C_{l,(h)} \bigr) = 0,    
\end{aligned}    
\end{equation}
where $E_{l,(h)} \in \mathbb{R}^{[m_1(r+1)]\times n}$ for $0 \leq l \leq r$, and
\[
    E_{l,(h)}=\begin{bmatrix} \mathbb{O}_{m_1 \times n}^{\times (l)} \\ \sqrt{-\frac{1}{r+1}R_{22}^\sharp\left( P^{(h)} \right)^{-1} }\left(N_{1,(h-1)}-R_{12}(P^{(h)})R_{22}(P^{(h)})^{-1}N_{2,(h-1)}\right) \\ \mathbb{O}_{m_1 \times n}^{\times (r-l)}\end{bmatrix}.
\]

By virtue of Proposition \ref{main_results:stabilizable_detectable}, the well-definedness of $Z^{(h)}$ is equivalent to verifying the mean-square stabilizability of the system
\[\left[A_{(h)},C_{1,(h)},\dotsb,C_{r,(h)};B_2,D_{1,2},\dotsb,D_{r,2}\right]\]
and the stochastic detectability of the system
\[\left[E_{0,(h)},E_{1,(h)},\dotsb,E_{r,(h)}; A_{(h)},C_{1,(h)},\dotsb,C_{r,(h)}\right].\]

To verify mean-square stabilizability, set $L^{(h)} = L-\hat{K}_2(P^{(h)})$, where $\hat{K}_2(P)$ is defined in \eqref{linearoperators_auxiliarymatrix}. The system $( A_{(h)} + B_2 L^{(h)}, C_{1,(h)} + D_{1,2} L^{(h)}, \dotsb, C_{r,(h)} + D_{r,2} L^{(h)} )$ is associated with $\mathcal{L}^*_{P^{(h)}}$. Since $\operatorname{Spec}\,\mathcal{L}^*_{P^{(h)}} \subset \mathbb{C}^{-}$, the required mean-square stabilizability follows immediately.

Define the auxiliary matrix
\[
\varTheta^{(h)} =\begin{bmatrix} \varTheta_{0}^{(h)} & \dotsb  & \varTheta_{r}^{(h)}\end{bmatrix}\in \mathbb{R}^{n\times [m_1(r+1)]},
\]
with
\[
\varTheta_{0}^{(h)}=-\left(B_{1}-B_{2}R_{22}(P^{(h)} )^{-1}R_{21}(P^{(h)} )\right)\sqrt{-(r+1)R_{22}^\sharp\left( P^{(h)} \right)^{-1} },
\]
\[
\varTheta_{l}^{(h)}=  -\left(D_{l1}-D_{l2}R_{22}(P^{(h)})^{-1}R_{21}(P^{(h)} )\right)\sqrt{-(r+1)R_{22}^\sharp\left( P^{(h)}  \right)^{-1} },\quad 1 \leq l \leq r.
\]
It can be confirmed that the system
\[(A_{(h)}+\varTheta_{0}^{(h)} E_{0,(h)}, C_{1,(h)}+\varTheta_{1}^{(h)} E_{1,(h)},\dotsb, C_{r,(h)}+\varTheta_{r}^{(h)} E_{r,(h)})\]
related to $\mathcal{L}^{(h-1,h)}$ is mean-square stable, which guarantees the stochastic detectability of 
\[
[E_{0,(h)},E_{1,(h)},\dotsb,E_{r,(h)}; A_{(h)},C_{1,(h)},\dotsb,C_{r,(h)}].
\]
Invoking Proposition \ref{main_results:stabilizable_detectable} again, ARE \eqref{proof_2} possesses a unique stabilizing solution $Z^{(h)}$ such that $R_{22}(P^{(h)}+Z^{(h)}) \succ 0$, and we have $\operatorname{Spec}\, \mathcal{L}^{(h,h+1)}\subset\mathbb{C}^{-}$.

In view of $\operatorname{Spec}\, \mathcal{L}^*_{P^{(h)}}\subset\mathbb{C}^{-}$, Lemma \ref{main_results:lemma} yields $\tilde{P}_{L} \succeq P^{(h)} + Z^{(h)}$ and $\operatorname{Spec}\, \mathcal{L}^*_{P^{(h+1)}}\subset\mathbb{C}^{-}$. The condition $R_{22}(P^{(h)} + Z^{(h)}) \succ 0$ and $R_{11}(P^{(h)} + Z^{(h)}) \preceq R_{11}(\tilde{P}_{L})\prec 0$ implies $P^{(h)} + Z^{(h)} \in \mathrm{Dom}\,\mathcal{G}$. Substituting \eqref{proof_2} into $\mathcal{G}( P^{(h)} + Z^{(h)} )$ via Proposition \ref{function_g:proposition_3} gives the desired expression of $\mathcal{G}(P^{(h)}+Z^{(h)})$.

Thus $a_{h}$--$c_{h}$ hold true. By mathematical induction, all properties $a_k$--$c_k$ are valid for all nonnegative integers $k$.

\textbf{Step 3. Convergence analysis.}
The sequence $\{P^{(k)}\}_{k \ge 0}$ is monotonically non-decreasing and upper bounded, hence it is convergent. Moreover,
\[
\lim_{h\rightarrow\infty}\sum_{s=h+1}^{\infty}Z^{(s)}=0,\qquad \lim_{k\rightarrow\infty}Z^{(k)}=0.
\]
Denote $P^*=\lim\limits_{k\rightarrow\infty}P^{(k)}$. By Proposition \ref{function_g:proposition_3}, we have
\begin{equation*}
    \begin{aligned}
    \mathcal{G}(P^*)=\mathcal{G}\left(\lim_{k\rightarrow\infty}P^{(k)}\right)
    =&\sum_{s=h+1}^{\infty}Z^{(s)}A_{(h)} + \sum_{s=h+1}^{\infty}A_{(h)}^{\top} Z^{(s)} + \sum_{s=h+1}^{\infty}\sum_{l=1}^{r}C_{l,(h)}^{\top} Z^{(s)} C_{l,(h)}\\
    &- \left( \sum_{s=h+1}^{\infty}B_2^{\top} Z^{(s)} + \sum_{s=h+1}^{\infty}\sum_{l=1}^{r}D_{l,2}^{\top} Z^{(s)} C_{l,(h)} \right)^{\top} \left( R_{22} + \sum_{l=1}^{r}D_{l,2}^{\top} P^* D_{l,2} \right)^{-1} \\
    &\quad \times \left( \sum_{s=h+1}^{\infty}B_2^{\top} Z^{(s)} + \sum_{s=h+1}^{\infty}\sum_{l=1}^{r}D_{l,2}^{\top} Z^{(s)} C_{l,(h)} \right)+\mathcal{G}\left(\sum_{s=0}^{h}Z^{(s)}\right)
    \end{aligned}
\end{equation*}
for all $h\geq0$. Since $\lim\limits_{h\rightarrow\infty}\sum_{s=h+1}^{\infty}Z^{(s)}=0$ and
\begin{equation*}
    \begin{aligned}
\lim_{h\rightarrow\infty}\mathcal{G}\left(\sum_{s=0}^{h}Z^{(s)}\right)
    =&\lim_{h\rightarrow\infty}-\left(N_{1,(h)}-R_{12}(P^{(h+1)})R_{22}(P^{(h+1)})^{-1}N_{2,(h)}\right)^{\top}R^{\sharp}_{22}(P^{(h+1)})^{-1}\\
    &\times \left(N_{1,(h)}-R_{12}(P^{(h+1)})R_{22}(P^{(h+1)})^{-1}N_{2,(h)}\right)=0,
    \end{aligned}
\end{equation*}
we conclude that $\mathcal{G}(P^*)=0$ and $P^* \in \overline{\mathbb{S}}_{+}^n$ solves the stochastic GTARE \eqref{zslqsdg:gtare}.

Let $\tilde{P}$ be the stabilizing solution of  the stochastic GTARE \eqref{zslqsdg:gtare}. Define $\tilde{L}=\hat{K}_2(\tilde{P})$ with $\hat{K}_2(\tilde{P})$ given in \eqref{linearoperators_auxiliarymatrix}. By Proposition \ref{function_g:proposition_2}, the stochastic GTARE \eqref{zslqsdg:gtare} can be transformed into the standard ARE \eqref{gtare_lw} associated with $\tilde{L}$, and $\tilde{P}$ is exactly its stabilizing solution. Following the stochastic detectability argument in Step 2 and combining Lemma \ref{main_results:stochastically_detectable_lemma}, we deduce that $(A_{\tilde{L}},C_{1\tilde{L}},\dotsb,C_{r\tilde{L}})$ is mean-square stable, which implies $\tilde{L} \in \mathcal{A}$.

Using Proposition \ref{function_g:proposition_3} and decomposing $P^* = (P^* - \tilde{P}) + \tilde{P}$, we expand $\mathcal{G}(P^*)$ as
\begin{align*}
\mathcal{G}(P^*)
=&(P^*-\tilde{P})\left(A + B_1K_1(\tilde{P})+B_2K_2(\tilde{P})\right)+\left(A + B_1K_1(\tilde{P})+B_2K_2(\tilde{P})\right)^{\top}(P^*-\tilde{P})\\
&+\sum_{l=1}^{r}\left(C_{l} + D_{l,1}K_1(\tilde{P})+D_{l,2}K_2(\tilde{P})\right)^{\top}(P^*-\tilde{P})\left(C_{l} + D_{l,1}K_1(\tilde{P})+D_{l,2}K_2(\tilde{P})\right)\\
&-N_2(\tilde{P},P^*-\tilde{P})^{\top} R_{22}(P^*)^{-1} N_2(\tilde{P},P^*-\tilde{P})-M^{\top}R_{22}^{\sharp}(P^*)^{-1}M+\mathcal{G}(\tilde{P}),
\end{align*}
where $M =  N_1(\tilde{P},P^*-\tilde{P})-R_{12}(P^*)R_{22}(P^*)^{-1}N_2(\tilde{P},P^*-\tilde{P})$, and $N_1(\tilde{P},P^*-\tilde{P})$ and $N_2(\tilde{P},P^*-\tilde{P})$ are defined in \eqref{main_results:operators_k&n}. 

It is obvious that $-M^{\top}R_{22}^{\sharp}(P^*)^{-1}M \succeq 0$.
From Proposition \ref{main_results:stabilizable_detectable}, we have $P^* \succeq \tilde{P}$. Notice that $\tilde{P}_{\tilde{L}} = \tilde{P}$, which yields $\tilde{P} \preceq P^* \preceq \tilde{P}_{\tilde{L}} = \tilde{P}$. Consequently, $P^*= \tilde{P}$, i.e., the iterative sequence converges to the unique stabilizing solution of the stochastic GTARE \eqref{zslqsdg:gtare}.
\end{proof}

\section{Numerical Example}
\label{sec:example}

To verify the effectiveness of the proposed computational method, a specific numerical example is provided, where the system parameters are randomly generated as follows:

\begin{align*}
A &= \begin{bmatrix}
        0.951718 & -0.270753 & -0.399002 \\
        -0.402652 & -0.332890 & -0.211269 \\
        0.865505 & 0.242714 & 0.051286
\end{bmatrix} &
Q &= \begin{bmatrix}
        3.030930 & -1.047214 & -0.749947 \\
        -1.047214 & 2.194723 & 1.675449 \\
        -0.749947 & 1.675449 & 3.333918
\end{bmatrix} \\[6pt]
C_1 &= \begin{bmatrix}
        0.428394 & -0.231063 & 1.819921 \\
        1.636217 & -0.730599 & 1.212696 \\
        -0.187034 & -0.538025 & 1.193243
\end{bmatrix} &
C_2 &= \begin{bmatrix}
        0.160803 & -1.563720 & -1.122872 \\
        -0.855569 & -0.072125 & 0.710778 \\
        0.091879 & 0.884612 & 0.977433
\end{bmatrix} 
\end{align*}
\begin{align*}
B_1 &= \begin{bmatrix}
        0.325576 & -0.086447 & -0.227595 \\
        -0.281475 & 0.347989 & -0.397287 \\
        -0.079578 & -0.101469 & -0.095411
\end{bmatrix} &
B_2 &= \begin{bmatrix}
        1.102509 & 0.029449 & 0.441306 \\
        0.188708 & 1.055548 & 0.363029 \\
        0.257869 & 0.036254 & 1.116683
\end{bmatrix} \\[6pt]
D_{11} &= \begin{bmatrix}
        0.007102 & 0.008817 & 0.004927 \\
        0.006974 & 0.000950 & 0.001090 \\
        0.009301 & 0.004565 & 0.001537
\end{bmatrix} &
D_{12} &= \begin{bmatrix}
        0.005681 & 0.002687 & 0.000159 \\
        0.000953 & 0.009114 & 0.008537 \\
        0.007980 & 0.009017 & 0.006756
\end{bmatrix}\\[6pt]
D_{21} &= \begin{bmatrix}
        0.009843 & 0.001642 & 0.003074 \\
        0.002716 & 0.001324 & 0.004221 \\
        0.008972 & 0.003174 & 0.003310
\end{bmatrix} &
D_{22} &= \begin{bmatrix}
        0.000377 & 0.005941 & 0.009387 \\
        0.003089 & 0.006808 & 0.007371 \\
        0.005642 & 0.006324 & 0.007449
\end{bmatrix} 
\end{align*}
\begin{align*}
R_{11} &= \begin{bmatrix}
        -1.658649 & -0.888904 & -0.707647 \\
        -0.888904 & -1.960809 & -0.586128 \\
        -0.707647 & -0.586128 & -1.471209
\end{bmatrix} &
R_{12} &= \begin{bmatrix}
        0.319075 & 0.370440 & 0.509807 \\
        0.212195 & 0.132083 & 0.857426 \\
        0.313025 & 0.677299 & 0.000364
\end{bmatrix} \\[6pt]
R_{21} &= \begin{bmatrix}
        0.319075 & 0.212195 & 0.313025  \\
        0.370440 & 0.132083 & 0.677299 \\
        0.509807 & 0.857426 & 0.000364
\end{bmatrix} &
R_{22} &= \begin{bmatrix}
        1.078023 & 0.529765 & 0.411900 \\
        0.529765 & 0.908261 & 0.312832 \\
        0.411900 & 0.312832 & 1.135871
\end{bmatrix} \\[6pt]
S_1 &= \begin{bmatrix}
        0.779145 & 0.430019 & 0.484125 \\
        0.792295 & 0.346524 & 0.258780 \\
        0.064140 & 0.841675 & 0.116926
\end{bmatrix} & 
S_2 &= \begin{bmatrix}
        0.616144 & 0.362532 & 0.803394 \\
        0.351730 & 0.980676 & 0.676821 \\
        0.513458 & 0.405686 & 0.268156
\end{bmatrix}
\end{align*}

After 13 external iterations, the iterative sequence converges to the stabilizing solution of the problem within the specified error tolerance $\epsilon=1\times 10^{-12}$. The norm of the equation residual is $4.9409 \times 10^{-6}$. The specific values of the stabilizing solution to the stochastic GTARE \eqref{zslqsdg:gtare} (retained to five decimal places) are presented as follows:
\[
\begin{bmatrix}
4.49040 & 0.33056 & -1.00210 \\
0.33056 & 4.89533 & -0.78154 \\
-1.00210 & -0.78154 & 8.61632
\end{bmatrix}.
\]

\begin{figure}[htbp]
\centering
\includegraphics[width=0.99\textwidth]{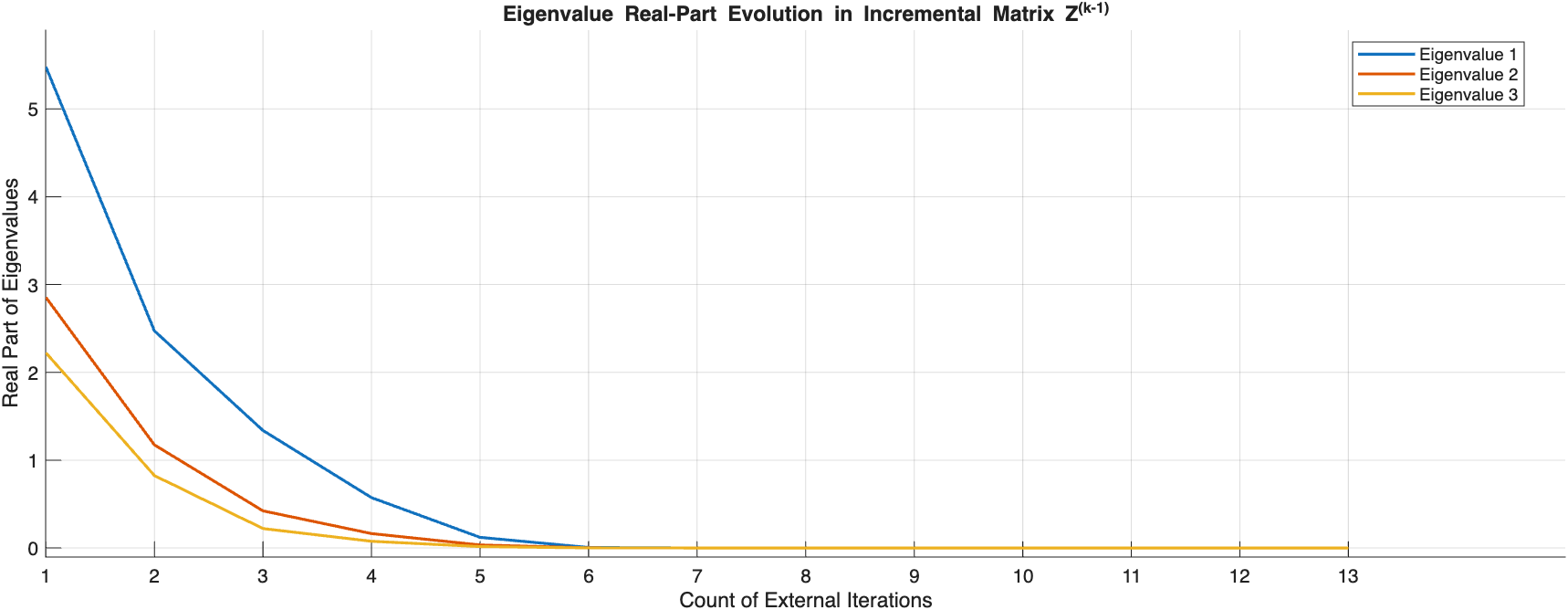}
\caption{The figure shows the three eigenvalues of $Z^{(k-1)}\,(k=1,\dotsb,13)$ evolving with external iterations.}
\label{fig:Figure_Z}
\end{figure}

\begin{figure}[htbp]
\centering
\includegraphics[width=0.99\textwidth]{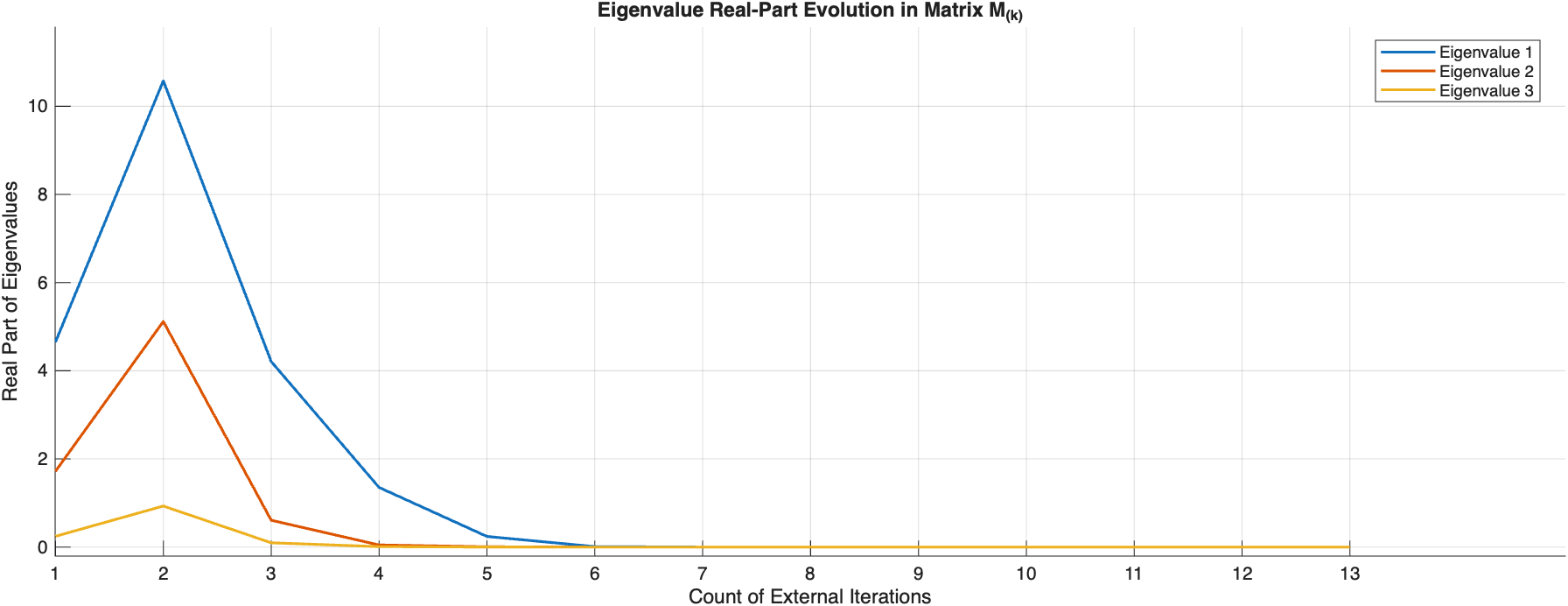}
\caption{The figure shows the three eigenvalues of $M_{(k-1)}(k=1,\dotsb,13)$ evolving with external iterations.}
\label{fig:Figure_M}
\end{figure}

\begin{figure}[htbp]
\centering
\includegraphics[width=0.99\textwidth]{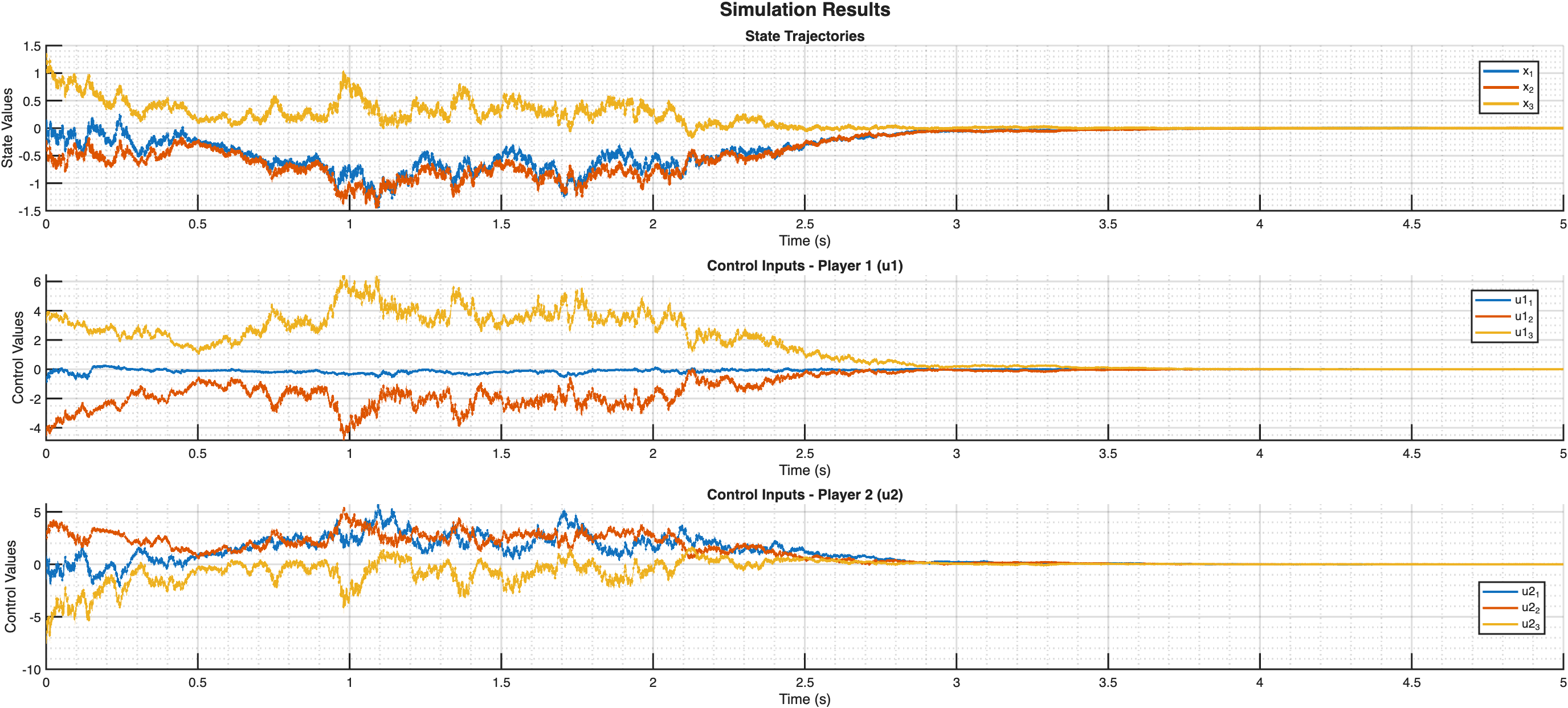}
\caption{This simulation diagram illustrates the dynamic evolution of the system state under the saddle-point equilibrium feedback strategies of the two players.}
\label{fig:Figure_S}
\end{figure}

Fig. \ref{fig:Figure_Z} and Fig. \ref{fig:Figure_M} illustrate the evolution of the eigenvalues of $Z^{(k-1)}$ and $M_{(k-1)}$ across external iterations $(k=1,\dotsb,13)$. Both figures reveal a consistent trend: the eigenvalues of $Z^{(k)}$ and $M_{(k)}$ decrease with iterations, approaching $0$ by the sixth iteration and continuing to decline thereafter. Notably, the eigenvalues of $M_{(1)}$ undergo a drastic change, which reflects the first correction step after the initial decoupling.
These observations admit a clear game-theoretic interpretation: when the set $\mathcal{A}$ is non-empty, the initial feedback gain $K_2(0)$ (associated with Player~2, the minimizer) already provides mean-square stabilization, while the corresponding gain $K_1(0)$ for Player~1 (the maximizer) is not yet at its best response. At each outer iteration, Player~2's feedback is updated not only to maintain system stability but also to rebalance against Player~1's dissatisfaction (measured by the residual $M_{(k)}$) from the previous step; the semidefinite correction $Z^{(k)}$ then captures Player~1's incremental improvement. Through this nested adjustment, the residual $\mathcal{G}(P^{(k)})$ monotonically decreases and the iterates approach the saddle-point equilibrium.

Fig. \ref{fig:Figure_S} illustrates the dynamic evolution of the system’s three state variables ($x_1$, $x_2$, $x_3$) starting from a randomly generated initial state. Initially, the system displays pronounced unstable dynamics, evidenced by the rapid and intense transients of all states—with $x_3$ exhibiting a particularly large initial deviation. This unstable trend is quickly suppressed, however, and the states converge asymptotically to the origin.  
The control input diagrams in the middle and lower panels explain this phenomenon: 
the controls applied by the two players take large values in opposite directions. Within the set of stabilizing feedback gains $\mathcal{K}$, the two players pursue mutually acceptable equilibrium control strategies to ensure the system’s equilibrium across all directions. While their control injection are antagonistic, their interaction yields a combined effect that counteracts the system’s inherent unstable dynamics, collectively driving the states toward equilibrium.  
Notably, despite the zero-sum game framework, the two parties demonstrate a cooperative tendency in maintaining system stability. These observations align with the theoretical findings presented in Remark 2.6.3 of \cite{Sun2020_book}.  

Together, these figures verify the effectiveness of our computational method and indirectly reveal the mechanism by which the method iteratively finds the equilibrium solution.

% \section*{Acknowledgments}

\clearpage
\printbibliography

\end{document}